\documentclass[10pt]{article}
\usepackage{amssymb,latexsym}
\usepackage{epsfig}
\usepackage{eufrak}
\usepackage{amsmath}
\usepackage{mathrsfs}
\usepackage{color}

\newtheorem{theorem}{Theorem}
\newtheorem{lemma}{Lemma}
\newtheorem{corollary}{Corollary}

\newcommand{\be}{\begin{equation}}
\newcommand{\ee}{\end{equation}}
\newcommand{\bee}{\begin{eqnarray*}}
\newcommand{\eee}{\end{eqnarray*}}
\newcommand{\bel}{\begin{eqnarray}}
\newcommand{\eel}{\end{eqnarray}}
\newcommand{\bec}{\begin{cases}}
\newcommand{\eec}{\end{cases}}
\newcommand{\bem}{\begin{bmatrix}}
\newcommand{\eem}{\end{bmatrix}}

\newcommand{\la}{\label}
\newcommand{\li}{\left}
\newcommand{\ri}{\right}

\newcommand{\DEF}{\stackrel{\mathrm{def}}{=}}

\newcommand{\lc}{\lceil}
\newcommand{\rc}{\rceil}

\newcommand{\ep}{\epsilon}
\newcommand{\vep}{\varepsilon}
\newcommand{\lm}{\lambda}

\newcommand{\Up}{\Upsilon}

\newcommand{\si}{\sigma}

\newcommand{\de}{\delta}

\newcommand{\vDe}{\varDelta}

\newcommand{\ga}{\gamma}

\newcommand{\se}{\theta}
\newcommand{\Se}{\Theta}

\newcommand{\ze}{\zeta}
\newcommand{\al}{\alpha}
\newcommand{\ba}{\beta}

\newcommand{\vro}{\varrho}
\newcommand{\ro}{\rho}

\newcommand{\om}{\omega}
\newcommand{\Om}{\Omega}

\newcommand{\f}{\frac}
\newcommand{\sq}{\sqrt}
\newcommand{\cd}{\cdots}

\newcommand{\qu}{\quad}
\newcommand{\qqu}{\qquad}
\newcommand{\fa}{\forall}

\newcommand{\mscr}{\mathscr}
\newcommand{\mcal}{\mathcal}

\newcommand{\bb}{\mathbb}

\newcommand{\mrm}{\mathrm}
\newcommand{\bs}{\boldsymbol}

\newcommand{\sh}{\slash}

\newcommand{\tx}{\text}

\newcommand{\iy}{\infty}

\newcommand{\bed}{\begin{description}}
\newcommand{\eed}{\end{description}}
\newcommand{\bei}{\begin{itemize}}
\newcommand{\eei}{\end{itemize}}
\newcommand{\ben}{\begin{enumerate}}
\newcommand{\een}{\end{enumerate}}
\newcommand{\bib}{\bibitem}
\newcommand{\beL}{\begin{lemma}}
\newcommand{\eeL}{\end{lemma}}
\newcommand{\beT}{\begin{theorem}}
\newcommand{\eeT}{\end{theorem}}

\newcommand{\beC}{\begin{corollary}}
\newcommand{\eeC}{\end{corollary}}

\newcommand{\sect}{\section}

\newcommand{\bpf}{\begin{pf}}
\newcommand{\epf}{\end{pf}}
\newcommand{\bsk}{\bigskip}

\newcommand{\pfbox}{\hfill\mbox{$\Box$}}

\newenvironment{pf}{\paragraph*{Proof{\rm.}}}{\pfbox\bigskip}

\begin{document}

\title{{\bf New Optional Stopping Theorems and Maximal Inequalities on Stochastic Processes} \thanks{The author had been previously working with Louisiana
State University at Baton Rouge, LA 70803, USA, and is now with Department of Electrical Engineering, Southern University and A\&M College,
Baton Rouge, LA 70813, USA; Email: chenxinjia@gmail.com. The main results of this paper have appeared in Proceedings of SPIE Conferences,
Baltimore, Maryland,  April 24-27, 2012.} }

\author{Xinjia Chen}

\date{First submitted in July 16, 2012}

\maketitle

\begin{abstract}

In this paper, we develop new optional stopping theorems  for scenarios where the stopping rules are defined by  bounded continuity regions.
Moreover, we establish a wide variety of inequalities on the supremums and infimums of functions of stochastic processes over the whole range of
time indexes.

\end{abstract}


\section{Introduction}

Martingale theory has been developed as a powerful tool for investigating stochastic processes.   One of the most important results of
martingale theory are Doob's optional stopping theorem \cite{Doob} and its variations.  These optional stopping theorems are relied on the
assumptions such as uniform integrability or integrable stopping times.  However, in many applications, the relevant stochastic process is not
uniformly integrable and the expectation of the stopping time is not necessarily finite.   Motivated by this situation, in this paper, we shall
develop new optional stopping theorems for scenarios where the uniform integrability of the stochastic process or the integrability of the
stoping time are not guaranteed, while the continuity region associated with the stopping rule is bounded. Based on the new optional stopping
theorems, we have established general maximal inequalities, which accommodate some classical inequalities such as, Bernstein's inequality
\cite{Bernstein}, Chernoff bounds \cite{Chernoff, Chernoffb}, Bennett's inequality \cite{Bennett}, Hoeffding-Azuma inequality \cite{Azuma,
Hoeffding} as special cases.

This paper is organized as follows. In Section 2, we present new optional stopping theorems.  In Section 3, we propose new maximal inequalities.
Section 3 is the conclusion.  All proofs are given in the Appendices.   The main results of this paper have appeared in \cite{Chen}.

Throughout this paper, we shall use the following notations. Let $\bb{R}$ denote the set of real numbers.  Let $\bb{R}^+$ denote the set of
non-negative real numbers. Let $\bb{Z}^+$ denote the set of non-negative integers.  Let $\bb{N}$ denote the set of positive integers.   Let
$(X_t)_{t \in \bb{T}}$ denote a stochastic process, where $\bb{T} \subseteq \bb{R}^+$ is the set of time values.  Specially, $(X_t)_{t \in
\bb{T}}$ is a continuous-time stochastic process if $\bb{T} = \bb{R}^+$; and $(X_t)_{t \in \bb{T}}$ is a discrete-time stochastic process if
$\bb{T} = \bb{Z}^+$. We assume that all stochastic processes are defined in probability space $(\Om, \mscr{F}, \Pr )$.  We also use $\bb{P}$ to
denote the probability measure $\Pr$.  For $t \in \bb{T}$, let $\mscr{F}_t$ denote the sub-$\si$-algebra generated by the collection of random
variables $\{ X_\tau: 0 \leq \tau \leq t, \; \tau \in \bb{T} \}$. The collection $(\mscr{F}_t)_{t \in \bb{T}}$ of sub-$\si$-algebras of
$\mscr{F}$ is called the natural filtration of $\mscr{F}$.    Let ``$A \vee B$'' denote the maximum of $A$ and $B$.  Let ``$A \wedge B$'' denote
the minimum of $A$ and $B$.  The other notations will be made clear as we proceed.

\sect{Optional Stopping Theorems}

In this section, we shall first develop some new optional stopping theorems on stochastic processes.  Consider a stochastic process $(X_t)_{t
\in \bb{R}^+}$ defined in the probability space $(\Om, \mscr{F}, \Pr )$. Let $\bs{\tau}$ be a stopping time taking values in $\bb{R}^+ \cup
\{\iy\}$. Define \[ X_{\bs{\tau}} = \lim_{t \to \iy} X_{\bs{\tau} \wedge t } \] if the limit exists. Clearly, for $\om \in \Om$, \[
X_{\bs{\tau}} (\om) = \bec X_{\bs{\tau} (\om)} (\om) & \tx{if $\bs{\tau} (\om) < \iy$};\\
\lim_{t \to \iy} X_t (\om) & \tx{if $\bs{\tau} (\om) = \iy$ and the limit exists}. \eec \] Let $\bs{\tau}_1$ and $\bs{\tau}_2$ be two stopping
times.   Since the stopping times can be $\iy$, we shall define the notion of  $\bs{\tau}_1 \leq \bs{\tau}_2$ as follows: \bee \{ \bs{\tau}_1
\leq \bs{\tau}_2 \} & = & \{\om \in \Om: \bs{\tau}_1 (\om) \leq \bs{\tau}_2 (\om), \; \bs{\tau}_1 (\om) \in \bb{R}^+, \; \bs{\tau}_2
(\om) \in \bb{R}^+ \}\\
&  &  \cup \{\om \in \Om:  \bs{\tau}_1 (\om) \in \bb{R}^+, \; \bs{\tau}_2 (\om) = \iy \} \cup \{\om \in \Om:  \bs{\tau}_1 (\om) = \iy, \;
\bs{\tau}_2 (\om) = \iy \}. \eee

Clearly, a discrete-time process $(X_k)_{k \in \bb{Z}^+}$ can be viewed as a right-continuous process $(X_t)_{t \in \bb{R}^+}$ with $X_t = X_k$
for $t \in [k, k+1), \; k \in \bb{Z}^+$.   Therefore, we shall consider the optional stopping problems in the general setting of continuous-time
processes.  However, in order to develop new optional stopping theorems for continuous-time processes, we first need to establish discrete-time
optional stopping theorems and then generalize them to continuous-time processes.  For a discrete-time process, we have the following general
results.

\beT \la{Themmain} Let $(X_k, \mscr{F}_k)_{k \in \bb{Z}^+ }$ be a discrete-time super-martingale.  Let $\bs{\tau}_1$ and $\bs{\tau}_2$ be two
stopping times such that $\bs{\tau}_1 \leq \bs{\tau}_2$ almost surely and that there exists a constant $C$ so that $\{ \bs{\tau}_2 > k \}
\subseteq \{ | X_k | < C \}$ for all $k \in \bb{Z}^+$.  Assume that \be \la{asmp888} \tx{ $X_{ \bs{\tau}_2}$ exist and  $\bb{E} [ | X_{
\bs{\tau}_2 } | ]$ is finite}. \ee Then, $\bb{E} [ X_{ \bs{\tau}_2  } \mid \mscr{F}_{ \bs{\tau}_1 } ] \leq X_{ \bs{\tau}_1 }$ and $\bb{E} [ X_{
\bs{\tau}_2 } ] \leq \bb{E} [ X_{ \bs{\tau}_1 } ]$ almost surely, with equality if $(X_k, \mscr{F}_k)_{k \in \bb{Z}^+}$ is a martingale.
Specially,  the assumption (\ref{asmp888}) is satisfied and the conclusion follows in the following cases:

(i) $(X_k, \mscr{F}_k)_{k \in \bb{Z}^+}$ is a super-martingale such that there exists a constant $\vDe$ so that $| X_{k+1} - X_k | < \vDe$
almost surely for all $k \in \bb{Z}^+$.

(ii) $(X_k, \mscr{F}_k)_{k \in \bb{Z}^+}$ is a non-negative super-martingale. \eeT

See Appendix \ref{Themmain_app} for a proof.

In particular, as an immediate application of Theorem \ref{Themmain}, we have the following result.

\beC \la{Them4} Let $(X_k)_{k \in \bb{Z}^+}$ be a discrete-time stochastic process such that $X_0, X_n - X_{n-1}, \; n = 1, 2, \cd$ are
independent random variables with zero means and that $\bb{E}[|X_0|] < \iy, \; \sup_{n > 0} |X_n - X_{n-1}| < \iy$ almost surely. Let
$\bs{\tau}_1$ and $\bs{\tau}_2$ be two stopping times such that $\bs{\tau}_1 \leq \bs{\tau}_2$ almost surely and that there exists a constant
$C$ so that $\{ \bs{\tau}_2 > n \} \subseteq \{ |X_n | < C \}$ for all $n \geq 0$. Then, $\bb{E} [ X_{ \bs{\tau}_2 } \mid \mscr{F}_{ \bs{\tau}_1
} ] = X_{ \bs{\tau}_1 }$ and $\bb{E} [ X_{ \bs{\tau}_2 } ] = \bb{E} [ X_{ \bs{\tau}_1 } ]$ almost surely.  \eeC

To investigate optional stopping problems for continuous-time processes, throughout the remainder of this paper, we shall define two collections
of sets of real numbers, denoted by $\{ \mscr{D}_t, \; t \in \bb{R}^+ \}$ and $\{ \mcal{D}_t, \; t \in \bb{R}^+ \}$,  such that the following
requirements are satisfied:

(i) $\mscr{D}_t \subseteq \mcal{D}_t$ for all $t \in \bb{R}^+$.

(ii) There exists a positive constant $C$ such that $\mcal{D}_t \subseteq [-C, C]$ for all $t \in \bb{R}^+$.

(iii) For any right-continuous function $g(t) : \bb{R}^+ \mapsto \bb{R}$, \bee &   & \lim_{n \to \iy}
\inf \{ t \in \bb{S}_n: g(t) \notin \mscr{D}_t \} = \inf \{ t \in \bb{R}^+: g(t) \notin \mscr{D}_t \},\\
&   & \lim_{n \to \iy} \inf \{ t \in \bb{S}_n: g(t) \notin \mcal{D}_t \} = \inf \{ t \in \bb{R}^+: g(t) \notin \mcal{D}_t \},  \eee where
$\bb{S}_n = \{ k 2^{-n}: k \in \bb{Z}^+ \}$ for $n \in \bb{N}$.  Note that the infimums can be $\iy$.

Clearly, if a stopping time $\bs{\tau}$ is defined such that $\{ \bs{\tau} > t \}$ implies $\{ X_t \in \mcal{D}_t \}$, then the region of $(X_t,
t)$ for continuing observing $(X_t)_{t \in \bb{R}^+}$ is bounded.  In this sense, a stopping rule with such a stopping time is called a stopping
rule with {\it bounded continuity region}.

In many areas of engineering and sciences,  it is a frequent problem to investigate a stochastic process with bounded rate of variation.  For
this purpose, the following result is useful.

\beT \la{Them3} Let $(X_t, \mscr{F}_t)_{t \in \bb{R}^+ }$ be a right-continuous super-martingale such that there exist constants $\de$ and $\vDe
> 0$ so that $| X_{t^\prime} - X_t | < \vDe$ almost surely provided that $| t^\prime - t | \leq \de$.
Let $\bs{\tau}_1 = \inf \{ t \in \bb{R}^+ : X_t \notin \mscr{D}_t \}$ and $\bs{\tau}_2 = \inf \{ t \in \bb{R}^+ : X_t \notin \mcal{D}_t \}$.
Then, $\bb{E} [ X_{ \bs{\tau}_2  } \mid \mscr{F}_{ \bs{\tau}_1 } ] \leq X_{ \bs{\tau}_1 }$ and $\bb{E} [ X_{ \bs{\tau}_2 } ] \leq \bb{E} [ X_{
\bs{\tau}_1 } ]$ almost surely, with equality if $(X_t, \mscr{F}_t)_{t \in \bb{R}^+ }$ is a martingale. \eeT

See Appendix \ref{Them3_app} for a proof.

\beT \la{Them1} Let $(X_t, \mscr{F}_t)_{t \in \bb{R}^+ }$ be a right-continuous, non-negative super-martingale.  Let $\bs{\tau}_1 = \inf \{ t
\in \bb{R}^+ : X_t \notin \mscr{D}_t \}$ and $\bs{\tau}_2 = \inf \{ t \in \bb{R}^+ : X_t \notin \mcal{D}_t \}$.  Then, $\bb{E} [ X_{ \bs{\tau}_2
} \mid \mscr{F}_{ \bs{\tau}_1 } ] \leq X_{ \bs{\tau}_1 }$ and $\bb{E} [ X_{ \bs{\tau}_2 } ] \leq \bb{E} [ X_{ \bs{\tau}_1 } ]$ almost surely,
with equality if $(X_t, \mscr{F}_t)_{t \in \bb{R}^+ }$ is a martingale. \eeT

See Appendix \ref{Them1_app} for a proof.

 It should be noted that Theorem \ref{Them1} can be readily generalized to a right-continuous super-martingale which is bounded from below by a constant.
 Making use of Theorem \ref{Them1}, we have established Corollary \ref{Coro2} as follows.

\beC \la{Coro2} Let $\mcal{V}_t$ be a right-continuous function of $t \geq 0$.  Let $(X_t)_{t \in \bb{R}^+ }$ be a right-continuous stochastic
process such that \be \la{rand} \bb{E} [ \exp( s X_0 ) ] \leq \exp ( \varphi(s) \mcal{V}_0 ), \qqu \bb{E} [ \exp( s (X_t - X_\tau) ) \mid
\mscr{F}_\tau ] \leq \exp ( \varphi(s) ( \mcal{V}_t - \mcal{V}_\tau ) ) \ee almost surely for arbitrary $t \geq \tau \geq 0$ and $s \in (-a,
b)$, where $a$ and $b$ are positive numbers or infinity, and $\varphi(s)$ is a function of $s \in (- a, b)$. Define $Y_t = \exp ( s X_t -
\varphi(s) \mcal{V}_t )$ for $s \in (-a, b)$.  Let $\bs{\tau}_1 = \inf \{ t \in \bb{R}^+: Y_t \notin \mscr{D}_t \}$ and $\bs{\tau}_2 = \inf \{ t
\in \bb{R}^+: Y_t \notin \mcal{D}_t \}$. Then, $\bb{E} [ Y_{ \bs{\tau}_2 } \mid \mscr{F}_{ \bs{\tau}_1 } ] \leq Y_{ \bs{\tau}_1 }$ and $\bb{E} [
Y_{ \bs{\tau}_2 } ] \leq \bb{E} [ Y_{ \bs{\tau}_1 } ]$ almost surely, with equality if (\ref{rand}) holds  with equality almost surely for
arbitrary $t \geq \tau \geq 0$ and $s \in (-a, b)$. \eeC

\section{Maximal Inequalities}

By virtue of the above optional stopping theorems, we shall establish some general maximal inequalities.  With regard to a uniformly integrable
(UI) martingale process, we have discovered the following fact.

\beT \la{Them6UI} If $(X_t, \mscr{F}_t)_{t \in \bb{R}^+}$ is a right-continuous uniformly integrable martingale which converges almost surely to
a constant $c$, then $X_t$ is equal to the constant $c$ for all $t \geq 0$ almost surely. \eeT

See Appendix \ref{Them6UI_app} for a proof.

Theorem \ref{Them6UI} implies that a right-continuous non-constant UI martingale never converges to a constant. For a super-martingale
converging to a constant, we have the following results.

\beT \la{Them5} Let $(X_t, \mscr{F}_t)_{t \in \bb{R}^+ }$ be a right-continuous, non-negative super-martingale which converges almost surely to
a constant $c$. Then, $\Pr \li \{  \sup_{t \geq 0} X_t \geq \ga \ri \} \leq \f{ \bb{E} [  X_0 ] - c }{\ga - c}$ for any $\ga
> c$.  Specially, $\Pr \li \{  \sup_{t \geq 0} X_t \geq \ga \ri \}$ is equal to $\f{ \bb{E} [  X_0 ] - c }{\ga - c}$ and $1$ in accordance
with $\ga > \bb{E} [  X_0 ]$ and $\ga \leq \bb{E} [  X_0 ]$ under additional assumption that $(X_t, \mscr{F}_t)_{t \in \bb{R}^+ }$ is a
continuous martingale. \eeT

See Appendix \ref{Them5_app} for a proof.

By virtue of Corollary  \ref{Coro2} and Markov's inequality, the following Corollary \ref{CThem7} can be established.

\beC \la{CThem7} Let $\mcal{V}_t$ be a right-continuous function of $t \geq 0$.  Let $(X_t)_{t \in \bb{R}^+}$ be a right-continuous stochastic
process such that $\bb{E} [ \exp( s X_0 ) ] \leq \exp ( \varphi(s) \mcal{V}_0 )$ and $\bb{E} [ \exp( s (X_t - X_\tau) ) \mid \mscr{F}_\tau ]
\leq \exp ( \varphi(s) ( \mcal{V}_t - \mcal{V}_\tau ) )$ almost surely for arbitrary $t \geq \tau \geq 0$ and $s \in (-a, b)$, where $a$ and $b$
are positive numbers or infinity, and $\varphi(s)$ is a function of $s \in (- a, b)$. Define $Y_t = \exp ( s X_t - \varphi(s) \mcal{V}_t  )$ for
$s \in (-a, b)$. Then,  $\Pr \li \{ \sup_{t \geq 0} Y_t \geq \ga \ri \} \leq  \f{1}{\ga}$ for $\ga > 0$.  \eeC

As a direct consequence of Theorem \ref{Them5}, we have shown the following Corollary \ref{expexact}.

\beC \la{expexact} Let $\mcal{V}_t$ be a non-negative, continuous function of $t \geq 0$ such that the limit inferior of $\mcal{V}_{n+1} -
\mcal{V}_n$ with respect to $n \in \bb{N}$ is positive. Let $(X_t)_{t \in \bb{R}^+}$ be a continuous stochastic process such that $\bb{E} [
\exp( s X_0 ) ] = \exp ( \varphi(s) \mcal{V}_0 )$ and $\bb{E} [ \exp( s (X_t - X_\tau) ) \mid \mscr{F}_\tau ] = \exp ( \varphi(s) ( \mcal{V}_t -
\mcal{V}_\tau ) )$ almost surely for arbitrary $t \geq \tau \geq 0$ and $s \in (-a, b)$, where $a$ and $b$ are positive numbers or infinity, and
$\varphi(s)$ is a function of $s \in (- a, b)$. Define $Y_t = \exp ( s X_t - \varphi(s) \mcal{V}_t  )$ for $s \in (-a, 0) \cup (0, b)$. Then,
$\Pr \li \{ \sup_{t \geq 0} Y_t \geq \ga   \ri \} =  \f{1}{\ga}$ for $\ga \geq 1$.  \eeC

See Appendix \ref{expexact_app} for a proof.

Making use of Corollary \ref{CThem7}, we have developed the following results concerning stochastic processes.

\beT  \la{Them9} Let $\mcal{V}_t$ be a non-negative, right-continuous function of $t \in [0, \iy )$. Let $(X_t)_{t \in \bb{R}^+ }$ be a
right-continuous stochastic process such that $\bb{E} [ \exp ( s (X_{t^\prime} - X_t)  ) \mid \mscr{F}_t ] \leq \exp ( (\mcal{V}_{t^\prime} -
\mcal{V}_t) \varphi(s) )$ almost surely for arbitrary $t^\prime \geq t \geq 0$ and $s \in (-a, b)$, where $a$ and $b$ are positive numbers or
infinity, and $\varphi(s)$ is a non-negative function of $s \in (- a, b)$.  Let $\tau \geq 0, \; \ga > 0, \; \eta \geq 0$. Define $\mscr{A} = \{
s \in (0, a): \varphi(- s) \leq \ga s \}$ and $\mscr{B} = \{ s \in (0, b): \varphi(s) \leq \ga s \} $. Then, {\small \bel & & \Pr \li \{ \inf_{t
> 0} \li [ X_t - X_0 + \ga \mcal{V}_\tau + \f{ \varphi(- s) }{s} (\mcal{V}_t - \mcal{V}_\tau) \ri ] \leq  0 \ri \}
\leq \li [ \exp \li ( \varphi(- s) - \ga s \ri ) \ri ]^{\mcal{V}_\tau} \qqu \fa s \in (0, a), \qqu \qqu \qqu \la{gena12}\\
&  & \Pr \li \{ \sup_{t > 0 } \li [ X_t - X_0 - \ga \mcal{V}_\tau  - \f{ \varphi(s) }{s} (\mcal{V}_t - \mcal{V}_\tau) \ri ] \geq  0 \ri \} \leq
\li [ \exp \li ( \varphi(s) - \ga s \ri ) \ri ]^{\mcal{V}_\tau} \qqu \fa s \in (0, b), \la{genb12} \\
& & \Pr \li \{ \inf_{t > 0} \li [ X_t - X_0  + \ga (\mcal{V}_\tau \vee \mcal{V}_t) \ri ] \leq 0 \ri \} \leq \inf_{ s \in (0, a) } \li [ \exp
\li ( \varphi(-s) - \ga s \ri ) \ri ]^{\mcal{V}_\tau}, \la{discovera99812}\\
&  & \Pr \li \{ \sup_{t > 0} \li [ X_t - X_0 - \ga (\mcal{V}_\tau \vee \mcal{V}_t) \ri ] \geq 0 \ri \} \leq \inf_{ s \in (0, b) } \li [  \exp
\li ( \varphi(s) - \ga s \ri ) \ri ]^{\mcal{V}_\tau} \la{discoverb89912} \qqu \tx{and moreover,} \\
 & & \Pr \li \{ \inf_{ t > 0 } \; ( X_t - X_0 + \eta +
\ga \mcal{V}_t ) \leq 0 \ri \} \leq \inf_{ s \in \mscr{A} } e^{- \eta s},  \la{discoverb9899add} \\
&  & \Pr \li \{ \sup_{ t
> 0 } \; ( X_t - X_0 - \eta  - \ga \mcal{V}_t ) \geq 0 \ri \} \leq \inf_{ s \in \mscr{B} }  e^{- \eta s},  \la{discoverb989912}\\
& & \Pr \li \{ \inf_{ t > 0 } \li [ X_t - X_0 + \eta + \ga (\mcal{V}_\tau \vee \mcal{V}_t) \ri ] \leq 0 \ri \} \leq \inf_{ s \in \mscr{A} } e^{-
\eta s} \li [ \exp \li ( \varphi(-s) - \ga s \ri ) \ri ]^{\mcal{V}_\tau}, \la{discovera88312}\\
&  & \Pr \li \{ \sup_{ t > 0 } \li [ X_t - X_0 - \eta  - \ga (\mcal{V}_\tau \vee \mcal{V}_t) \ri ] \geq 0 \ri \} \leq \inf_{ s \in \mscr{B}  }
e^{- \eta s} \li [ \exp \li ( \varphi(s) - \ga s  \ri ) \ri ]^{\mcal{V}_\tau} \la{discoverb9812} \eel} provided that $\mscr{A}$ and $\mscr{B}$
are nonempty respectively.  In particular, under the above assumptions on $X_t, \; \mcal{V}_t$ and $\varphi(s)$, the following statements hold
true:

\noindent (I): If $\varphi(s)$ is a continuous function smaller than $\ga |s|$ at a neighborhood of $0$, then {\small \bel &  & \Pr  \{ \inf_{t
> 0} [ X_t - X_0 + \ga \mcal{V}_\tau  + \al(\ga) (\mcal{V}_t - \mcal{V}_\tau)  ] \leq 0 \} \leq \inf_{s \in (0, a) }  [ \exp ( \varphi(- s) -
\ga s ) ]^{\mcal{V}_\tau}, \la{discovera12}\\
&  & \Pr \{ \sup_{t > 0} [ X_t - X_0 - \ga \mcal{V}_\tau - \ba (\ga) (\mcal{V}_t - \mcal{V}_\tau) ] \geq 0 \} \leq \inf_{s \in (0, b) } [ \exp (
\varphi(s) - \ga s ) ]^{\mcal{V}_\tau}, \la{discoverb12} \eel} where $\al (\ga)$ and $\ba (\ga)$ are functions of $\ga$ defined as follows: $\al
(\ga)$ is equal to $\f{\varphi(- s^\star)}{s^\star}$ if $\inf_{s \in (0, a) } \li [ \varphi(- s) - \ga s \ri ]$ is attained at $s^\star \in (0,
a)$ and otherwise equal to $\lim_{s \uparrow a} \f{\varphi(- s)}{s}$; $\ba (\ga)$ is equal to $\f{\varphi(s^*)}{s^*}$ if $\inf_{s \in (0, b) }
\li [ \varphi(s) - \ga s \ri ]$ is attained at $s^* \in (0, b)$ and otherwise equal to $\lim_{s \uparrow b} \f{\varphi(s)}{s}$. Moreover, $0 <
\al (\ga) < \ga$ and $0 < \ba (\ga) < \ga$.

\noindent (II): If $\f{ \varphi(s) }{|s|}$ is monotonically increasing with respect to $|s| > 0$, then {\small \be \la{Julyb} \Pr\{ \inf_{t > 0}
[ X_t - X_0 + \eta + \ga (\mcal{V}_\tau \vee \mcal{V}_t) ] \leq 0 \} \leq \inf_{ s \in (0, a^\star) }  e^{- \eta s} [  \exp ( \varphi(-s) - \ga
s ) ]^{\mcal{V}_\tau} \ee } and {\small \be \la{Julya} \Pr \{ \sup_{t > 0} [ X_t - X_0 - \eta  - \ga (\mcal{V}_\tau \vee \mcal{V}_t) ] \geq 0 \}
\leq \inf_{ s \in (0, b^\star) }  e^{- \eta s} [ \exp \li ( \varphi(s) - \ga s \ri )  ]^{\mcal{V}_\tau}, \ee}  where $a^\star$ and $b^\star$ are
defined as follows: $a^\star$ is equal to $a$ if $\lim_{s \uparrow a} \f{\varphi(- s)}{s} \leq \ga$ and otherwise equal to  $s \in (0, a)$ such
that $\f{\varphi(- s)}{s} = \ga$;  $b^\star$ is equal to $b$ if $\lim_{s \uparrow b} \f{\varphi(s)}{s} \leq \ga$ and otherwise equal to  $s \in
(0, b)$ such that $\f{\varphi(s)}{s} = \ga$.

\eeT

See Appendix \ref{Them9_app} for a proof.

An important application of Theorem \ref{Them9} is illustrated as follows.  Let $Y$ be a random variable with mean $\mu$. Define $X_0 = 0$ and
$X_n = \sum_{i=1}^n (Y_i - \mu)$ for $n \in \bb{N}$, where $Y_1, Y_2, \cd$ are i.i.d. samples of $Y$. Define a right-continuous stochastic
process $(X_t)_{t \in \bb{R}^+}$ such that $X_t = X_n$ for $t \in [n, n+1), \; n = 0, 1, 2, \cd$. Define a right-continuous function
$\mcal{V}_t$ of $t \in \bb{R}^+$ such that $\mcal{V}_t = n$ for $t \in [n, n+1), \; n = 0, 1, 2, \cd$.  If there is a convex function
$\varphi(s)$ such that $\ln \bb{E} \li [ \exp \li ( s (Y-\mu) \ri ) \ri ] \leq \varphi(s)$ for $s \in (-a, b)$ and that $\varphi(0) = 0$, then
we can apply Theorem \ref{Them9} to develop maximal inequalities for $(X_t)_{t \in \bb{R}^+}$, which immediately lead to maximal inequalities
for $(X_n)_{n \in \bb{N}}$. The function $\varphi(s)$ of the desired properties can be found for some particular cases as follows:

\noindent (i) If the moment generating function of $Y$ exists, then $\varphi(s)$ can be taken as $\ln \bb{E} \li [ \exp \li ( s Y \ri ) \ri ] -
\mu s$.

\noindent (ii) If $Y$ is a random variable such that $\bb{E} [ Y ] = 0, \; \bb{E} [ Y^2] = \si^2$ and $Y \leq b$, then by Bennett's inequality
\cite{Bennett}, the function $\varphi(s)$ can be taken as {\small $\ln [  \f{b^2}{b^2 + \si^2} \exp ( - \f{\si^2}{b} s ) + \f{\si^2}{b^2 +
\si^2} \exp (b s) ]$}.

\noindent (iii) If $Y$ is a random variable bounded in interval $[0, 1]$ almost surely,  then by Hoeffding's inequality \cite{Hoeffding}, the
function $\varphi(s)$ can be taken as $\ln ( 1 - \mu + \mu e^s ) - \mu s$.

\noindent (iv) If $Y$ is  a random variable uniformly distributed over $[- \f{1}{2}, \f{1}{2}]$, then we can show that $\bb{E} [ e^{s Y}] \leq
\exp ( \f{s^2}{24} )$ for all real number $s$.  Hence, the function $\varphi(s)$ can be taken as $\f{s^2}{24}$.  See Appendix \ref{app_uniform}
for the development of the bound for the moment generating function $\bb{E} [ e^{s Y}]$.

Applying Corollary \ref{expexact} to a continuous stochastic process, we have the following  results.

\beT  \la{Them10} Let $\mcal{V}_t$ be a non-negative, continuous function of $t \in [0, \iy )$ such that the limit inferior of $\mcal{V}_{n+1} -
\mcal{V}_n$ with respect to $n \in \bb{N}$ is positive. Let $(X_t)_{t \in \bb{R}^+}$ be a continuous stochastic process such that $\bb{E} [ \exp
( s (X_{t^\prime} - X_t) ) \mid \mscr{F}_t ] = \exp ( (\mcal{V}_{t^\prime} - \mcal{V}_t) \varphi(s) )$ almost surely for arbitrary $t^\prime
\geq t \geq 0$ and $s \in (-a, b)$, where $a$ and $b$ are positive numbers or infinity, and $\varphi(s)$ is a non-negative function of $s \in (-
a, b)$.  Let $\tau > 0, \; \ga > 0$ and $\eta > 0$. Then, {\small $\Pr \{ \inf_{t > 0} [ X_t - X_0 + \ga \mcal{V}_\tau  + \f{ \varphi(- s) }{s}
(\mcal{V}_t - \mcal{V}_\tau) ] \leq  0 \} = [ \exp ( \varphi(- s) - \ga s ) ]^{\mcal{V}_\tau}$} for any $s \in (0, a)$ and
 {\small $\Pr \{ \sup_{t > 0 } [ X_t - X_0 - \ga \mcal{V}_\tau  - \f{ \varphi(s) }{s} (\mcal{V}_t - \mcal{V}_\tau) ] \geq  0 \} =
 [ \exp ( \varphi(s) - \ga s ) ]^{\mcal{V}_\tau}$} for any $s \in (0, b)$.  In particular, under the above assumptions on $X_t, \;
\mcal{V}_t$ and $\varphi(s)$, the following statements hold true:

\noindent (I): If there exists $s^\star \in (0, a)$ such that $\varphi(- s^\star) = \ga s^\star$, then {\small $\Pr \li \{ \inf_{ t
> 0 } \; (  X_t - X_0 + \eta + \ga \mcal{V}_t ) \leq 0 \ri \}  = e^{- \eta s^\star}$}.

\noindent (II): If there exists $s^* \in (0, b)$ such that $\varphi(s^*) = \ga s^*$, then {\small $\Pr \li \{ \sup_{ t > 0 } \; ( X_t - X_0 -
\eta - \ga \mcal{V}_t ) \geq 0 \ri \} = e^{- \eta s^*}$}.

\noindent (III): If $\varphi(s)$ is a continuous function smaller than $\ga |s|$ at a neighborhood of $0$, then {\small $\Pr \{ \inf_{t > 0} [
X_t - X_0 + \ga \mcal{V}_\tau  + \al(\ga) (\mcal{V}_t - \mcal{V}_\tau) ] \leq 0 \} = \inf_{s \in (0, a) } [ \exp ( \varphi(- s) - \ga s )
]^{\mcal{V}_\tau}$} and {\small $\Pr \{ \sup_{t > 0} [ X_t - X_0 - \ga \mcal{V}_\tau - \ba (\ga) (\mcal{V}_t - \mcal{V}_\tau) ] \geq 0 \} =
\inf_{s \in (0, b) } [ \exp ( \varphi(s) - \ga s ) ]^{\mcal{V}_\tau}$},  where $\al (\ga)$ and $\ba (\ga)$ are functions of $\ga$ defined as
follows: $\al (\ga)$ is equal to $\f{\varphi(- s^\star)}{s^\star}$ if $\inf_{s \in (0, a) } \li [ \varphi(- s) - \ga s \ri ]$ is attained at
$s^\star \in (0, a)$ and otherwise equal to $\lim_{s \uparrow a} \f{\varphi(- s)}{s}$; $\ba (\ga)$ is equal to $\f{\varphi(s^*)}{s^*}$ if
$\inf_{s \in (0, b) } \li [ \varphi(s) - \ga s \ri ]$ is attained at $s^* \in (0, b)$ and otherwise equal to $\lim_{s \uparrow b}
\f{\varphi(s)}{s}$.  Moreover, $0 < \al (\ga) < \ga$ and $0 < \ba (\ga) < \ga$.

\eeT

Applying Theorem \ref{Them9} to i.i.d random variables with common probability density (or mass) function in an exponential family, we have
shown the following results, which generalize Chernoff bounds \cite{Chernoff}.

\beC  \la{expmax} Let $Y_1, Y_2, \cd$ be i.i.d. random samples of $Y$ which possesses a probability density (or mass) function $f_Y(y; \se) =
w(y) \exp ( u (\se) y - v(\se) )$ such that $\f{ d v( \se ) }{ d \se } = \se \f{ d u ( \se ) }{d \se}$ for $\se \in \Se$.  Define $ X_n
 = \sum_{i = 1}^n Y_i$ for $n \in \bb{N}$.
 Define $\mscr{M} (z, \se ) = \f{ \exp \li ( u (\se )   z -  v(\se) \ri ) } { \exp \li ( u (z ) z - v(z) \ri ) }$ and
 $\ro (z, \se, m, n) = m z + (n - m) \f{ v(z) - v(\se) }{ u (z) - u (\se) }$ for $z, \se \in \Se$ and $m, n \in \bb{N}$.
 Then, for all integer $m > 0$ and real number $\ga > 0$,
{\small  \bel &  &  \Pr \li \{ \sup_{n > 0} \li [ X_n - n \se  - \ga(n \vee m) \ri ] \geq  0 \ri \} \leq \Pr \li \{ \sup_{n > 0}
 \li [ X_n -  \ro (\se + \ga, \se, m, n) \ri ] \geq  0 \ri \} \leq  [ \mscr{M} (\se + \ga, \se ) ]^m  \qqu \la{ineqgood8}\\
&  & \Pr \li \{ \inf_{n > 0} \li [ X_n -  n \se  + \ga (n \vee m) \ri ] \leq 0 \ri \} \leq
 \Pr \li \{ \inf_{n > 0} \li [ X_n -  \ro (\se - \ga, \se, m, n) \ri ] \leq 0 \ri \} \leq [ \mscr{M} (\se - \ga, \se ) ]^m, \qqu \la{ineqgood9}
\eel} provided that $\se + \ga \in \Se$ and $\se - \ga \in \Se$ respectively.

\eeC

See Appendix \ref{expmax_app} for a proof.

By virtue of Theorem \ref{Them9}, we can generalize Hoeffding-Azuma's inequality  \cite{Azuma, Hoeffding} as follows.

\beC \la{CorAzuma} Let $(X_t)_{t \in \bb{R}^+}$ be a right-continuous stochastic process. Assume that there exist a right-continuous function
$\mcal{V}_t$ and a stochastic process $(Y_t)_{t \in \bb{R}^+}$ such that for all $t \geq 0$, $Y_t$ is measurable in $\mscr{F}_t$, and that $|X_t
- Y_\tau|^2 \leq \mcal{V}_t - \mcal{V}_\tau$ almost surely for arbitrary $t \geq \tau \geq 0$.  For $\ga > 0$ and $\tau > 0$ such that
$\mcal{V}_\tau > 0$, the following statements hold true:

\noindent (I) If $(X_t, \mscr{F}_t)_{t \in \bb{R}^+}$ is a super-martingale, then \be \la{hoe9a} \Pr \li \{   \sup_{t > 0} \li [  X_t - X_0 -
\f{\ga}{2} \li ( 1 + \f{\mcal{V}_t}{\mcal{V}_\tau} \ri ) \ri ] \geq 0 \ri \} \leq \exp \li ( - \f{\ga^2} { 2 \mcal{V}_\tau } \ri ). \ee

\noindent (II) If $(X_t, \mscr{F}_t)_{t \in \bb{R}^+}$ is a sub-martingale, then \be \la{hoe9b}  \Pr \li \{   \inf_{t > 0} \li [  X_t - X_0 +
\f{\ga}{2} \li ( 1 + \f{\mcal{V}_t}{\mcal{V}_\tau} \ri ) \ri ] \leq 0   \ri \} \leq \exp \li ( - \f{\ga^2} { 2 \mcal{V}_\tau } \ri ). \ee

\noindent (III) If $(X_t, \mscr{F}_t)_{t \in \bb{R}^+}$ is a martingale, then \be \la{hoe9c} \Pr \li \{   \sup_{t > 0} \li [  \li | X_t - X_0
\ri | - \f{\ga}{2} \li ( 1 + \f{\mcal{V}_t}{\mcal{V}_\tau} \ri ) \ri ] \geq 0 \ri \} \leq2 \exp \li ( - \f{\ga^2} { 2 \mcal{V}_\tau } \ri ). \ee

\eeC

See Appendix \ref{CorAzuma_app} for a proof.

With the help of Theorem \ref{Them9}, we have generalized  Bernstein's inequality \cite{Bernstein}, Bennett's  inequality \cite{Bennett} and
Chernoff bound \cite{Chernoffb} as follows.

\beC

\la{CBB}

Let $(X_k, \mscr{F}_k)_{k \in \bb{Z}^+ }$ be a martingale satisfying $\mrm{Var} (X_n \mid \mscr{F}_{n - 1} ) \DEF  \bb{E} [ ( X_n - \bb{E}[X_n
\mid \mscr{F}_{n - 1}] )^2 \mid \mscr{F}_{n - 1} ] \leq \si_n^2$ and $X_n - X_{n - 1} \leq a_n + b$ almost surely for $n \in \bb{N}$, where $b >
0$ and $a_n$ are deterministic numbers. Define $\mcal{V}_n = \sum_{i = 1}^n (\si_i^2 + a_i^2)$ for $n \in \bb{N}$.  Then,  {\small \bel  \Pr \li
\{ \sup_{n > 0} \li [ X_n - X_0 - \ga \mcal{V}_m - \li ( \f{\ga} {\ln (1 + b \ga)} - \f{1}{b} \ri )  (\mcal{V}_n - \mcal{V}_m) \ri ] \geq  0 \ri
\} \leq \li [ \exp \li ( \f{\ga}{b} - \f{(1 + b \ga) \ln (1
+ b \ga) }{b^2} \ri ) \ri ]^{ \mcal{V}_m},  \qqu  \la{beneet} &  & \\
\Pr  \li \{ \sup_{n > 0} \li [ X_n - X_0 - \f{\ga}{2} \li ( 1 + \f{\mcal{V}_n}{\mcal{V}_m}  \ri ) \ri ] \geq 0 \ri \} \leq \exp \li ( -
\f{\ga^2}{ 2 ( \mcal{V}_m + b \ga \sh 3) } \ri ) \qqu \qqu \qqu \qqu \qqu \qqu \qqu \qqu \qqu \la{Bert} &  & \eel} for all integer $m > 0$ and
real number $\ga
> 0$. Specially, if $b = 1$ and $a_n = 0$ for $n \in \bb{N}$, then \be \la{chert8} \Pr \li \{ \sup_{n > 0} \li [ X_n - X_0 - \f{\ga}{2} \li ( 1
+ \f{\mcal{V}_n}{\mcal{V}_m} \ri ) \ri ] \geq 0 \ri \} \leq \exp \li ( - \f{\ga^2}{ 4 \mcal{V}_m } \ri ) \ee for all integer $m
> 0$ and real number $\ga \in ( 0, \f{7}{2} \mcal{V}_m)$. \eeC

See Appendix \ref{CBB_app} for a proof.

Applying Theorem \ref{Them9} to a Poisson process, we have obtained the following results.

\beC \la{Posmax} Let $X_t$ be the number of arrivals in time interval $[0, t]$ for a Poisson process with an arrival rate $\lm > 0$.  Then,
{\small \bee &  & \Pr \li \{ \sup_{t > 0} \li [ X_t - (\lm + \ga) \tau - \f{\ga (t - \tau) }{ \ln (1 + \f{\ga}{\lm } ) }  \ri ] \geq 0 \ri \}
\leq \li  [ \li ( \f{\lm}{\lm + \ga} \ri )^{\lm + \ga} e^\ga \ri ]^\tau,\\
&  & \Pr \li \{  \inf_{t > 0} \li [ X_t - (\lm - \ga) \tau + \f{ \ga (t - \tau) }{ \ln (1 - \f{\ga}{\lm } ) }   \ri ] \leq 0 \ri \} \leq \li [
\li ( \f{\lm}{ \lm - \ga } \ri  )^{\lm - \ga} e^{-\ga} \ri ]^\tau \eee} for any $\tau > 0$ and $\ga > 0$. \eeC

See Appendix \ref{Posmax_app} for a proof.

\section{Conclusion}

In this paper, we have developed some new optional stopping theorems on martingale processes which require no assumption of uniform
integrability and integrable stopping times.  Making use of bounds of moment generating functions of increments of stochastic processes, we have
established a wide class of maximal inequalities on stochastic processes, which includes classical results such as Chernoff bounds,
Hoeffding-Azuma inequalities as special cases.

\appendix

\section{Proof of Theorem \ref{Themmain}} \la{Themmain_app}

Throughout the proof of the theorem, let $A \in \mscr{F}_{ \bs{\tau}_1 }$ and $B = A \cap  \{ \bs{\tau}_1 = n \} \in \mscr{F}_n$.
\subsection{Proof of Discrete-Time Optional Stopping Theorem under Assumption (\ref{asmp888}) }

In this section of Appendix \ref{Themmain_app},  we shall show the discrete-time optional stopping theorem under assumption (\ref{asmp888}).
More formally, we want to prove the following result:

{\it Let $(X_k, \mscr{F}_k)_{k \in \bb{Z}^+ }$ be a super-martingale.  Let $\bs{\tau}_1$ and $\bs{\tau}_2$ be two stopping times such that
$\bs{\tau}_1 \leq \bs{\tau}_2$ almost surely and that there exists a constant $C$ so that $\{ \bs{\tau}_2 > k \} \subseteq \{ | X_k | < C \}$
for all $k \in \bb{Z}^+$.  Assume that $X_{ \bs{\tau}_2}$ exist and  $\bb{E} [ | X_{ \bs{\tau}_2 } | ]$ is finite.  Then, $\bb{E} [ X_{
\bs{\tau}_2  } \mid \mscr{F}_{ \bs{\tau}_1 } ] \leq X_{ \bs{\tau}_1 }$ and $\bb{E} [ X_{ \bs{\tau}_2 } ] \leq \bb{E} [ X_{ \bs{\tau}_1 } ]$
almost surely, with equality if $(X_k, \mscr{F}_k)_{k \in \bb{Z}^+}$ is a martingale.}

\bsk

The following result stated as Lemma \ref{lem1a} is due to Doob \cite{Doob}, which can be found in many text books of probability theory.

\beL \la{lem1a} For $i \geq n$, \be \la{induce} \int_{ B \cap \{ \bs{\tau}_2 \geq n \}}  X_n  \; d \bb{P} \geq \int_{ B \cap \{ n \leq
\bs{\tau}_2 \leq i \}} X_{\bs{\tau}_2} \; d \bb{P} + \int_{ B \cap \{ \bs{\tau}_2 > i \}}  X_i  \; d \bb{P}, \ee  with equality if $X_t$ is a
martingale. \eeL

\beL \la{recall} $\Pr \li \{ \lim_{i \to \iy} \bb{I}_{ B \cap \{ n \leq \bs{\tau}_2 \leq i \}}  =   \bb{I}_{ B \cap \{ n \leq \bs{\tau}_2 < \iy
\}} \ri \} = 1$.

\eeL

\bpf

To show this,  consider two cases. In the case that $\om \notin B \cap \{ n \leq \bs{\tau}_2 < \iy \}$, we have $\om \notin B \cap \{ n \leq
\bs{\tau}_2 \leq i \}$ for all $i \geq n$ and thus $\bb{I}_{ B \cap \{ n \leq \bs{\tau}_2 < \iy \}} = 0 = \lim_{i \to \iy} \bb{I}_{ B \cap \{ n
\leq \bs{\tau}_2 \leq i \}}$. In the case that $\om \in B \cap \{ n \leq \bs{\tau}_2 < \iy \}$, we have $\bb{I}_{ B \cap \{ n \leq \bs{\tau}_2 <
\iy \}} = 1$. Then, $\bb{I}_{ B \cap \{ n \leq \bs{\tau}_2 \leq i \}} = 1$ for $i \geq \bs{\tau}_2 (\om)$, which implies that $\lim_{i \to \iy}
\bb{I}_{ B \cap \{ n \leq \bs{\tau}_2 \leq i \}} = 1$.  This proves the lemma.

\epf

\beL \la{lem889} Assume that $\bb{E} [ | X_{\bs{\tau}_2} | ] < \iy$. Then, \bee \lim_{i \to \iy} \int_{ B \cap \{ n \leq \bs{\tau}_2 \leq i \}}
X_{\bs{\tau}_2} \; d \bb{P} = \int_{ B \cap \{ n \leq \bs{\tau}_2 < \iy \}} X_{\bs{\tau}_2} \; d \bb{P} \leq \int_{ B \cap \{ n \leq \bs{\tau}_2
< \iy \}} | X_{\bs{\tau}_2} | \; d \bb{P} < \iy. \eee \eeL

\bpf From the assumption that $\bb{E} [ | X_{\bs{\tau}_2} | ] < \iy$, we have
\[
\bb{E} \li [ X_{\bs{\tau}_2}  \bb{I}_{ B \cap \{ n \leq \bs{\tau}_2 < \iy \}} \ri ] \leq \bb{E} \li [ | X_{\bs{\tau}_2} | \;  \bb{I}_{ B \cap \{
n \leq \bs{\tau}_2 < \iy \}} \ri ] \leq \bb{E} [ | X_{\bs{\tau}_2} | ] < \iy.
\]
By virtue of Lemma \ref{recall}, we have that $X_{\bs{\tau}_2} \bb{I}_{ B \cap \{ n \leq \bs{\tau}_2 \leq i \}} \to X_{\bs{\tau}_2}  \bb{I}_{ B
\cap \{ n \leq \bs{\tau}_2 < \iy \}}$ almost surely as $i \to \iy$.  Since $| X_{\bs{\tau}_2} \bb{I}_{ B \cap \{ n \leq \bs{\tau}_2 \leq i \}} |
\leq  | X_{\bs{\tau}_2} | \; \bb{I}_{ B \cap \{ n \leq \bs{\tau}_2 < \iy \}}$,  the lemma follows from the dominated convergence theorem.

\epf

\beL \la{lem8899}
\[
\lim_{i \to \iy} \int_{ B \cap \{ \bs{\tau}_2 > i \}}  X_i  \; d \bb{P} = \int_{ B \cap \{ \bs{\tau}_2 = \iy \}}  X_{\bs{\tau}_2}  \; d \bb{P}
\leq \int_{ B \cap \{ \bs{\tau}_2 = \iy \}}  | X_{\bs{\tau}_2}  | \; d \bb{P} < \iy.
\]
\eeL

\bpf

First, we shall show that \be \la{claim81}
 \lim_{i \to \iy} \int_{ B \cap \{ i < \bs{\tau}_2 < \iy \}}  X_i  \; d \bb{P}  = 0.
\ee
By the assumption that there exists a constant $C$ so that $\{ \bs{\tau}_2 > t \} \subseteq \{ |X_t | < C \}$ for any $t \geq 0$, we have
\bee \lim_{i \to \iy} \li | \int_{ B \cap \{ i < \bs{\tau}_2 < \iy \}}  X_i  \; d \bb{P} \ri | & \leq & \lim_{i \to \iy} \int_{ B \cap \{
i < \bs{\tau}_2 < \iy \}}  | X_i| \; d \bb{P}\\
&  \leq  & \lim_{i \to \iy} \int_{ B \cap \{ i < \bs{\tau}_2 < \iy \}} C \; d \bb{P} \leq \lim_{i \to \iy} \int_{ \{ i < \bs{\tau}_2 < \iy \}} C
\; d \bb{P} \\
& = & C \lim_{i \to \iy} \Pr \{ i < \bs{\tau}_2 < \iy \} = 0,  \eee which implies (\ref{claim81}).

Next, we shall show that \be \la{claim82} \lim_{i \to \iy} \int_{ B \cap \{ \bs{\tau}_2 = \iy \}}  X_i  \; d \bb{P} = \int_{ B \cap \{
\bs{\tau}_2 = \iy \}}  X_{\bs{\tau}_2}  \; d \bb{P}. \ee  By the assumption that there exists a constant $C$ so that $\{ \bs{\tau}_2 > t \}
\subseteq \{ |X_t | < C \}$ for any $t \geq 0$, we have $| X_i \; \bb{I}_{ B \cap \{ \bs{\tau}_2 = \iy \}} | \leq C$ for all $i \geq n$.  By the
definition of $X_{\bs{\tau}_2}$, we have that $X_i \; \bb{I}_{ B \cap \{ \bs{\tau}_2 = \iy \}} \to X_{\bs{\tau}_2} \; \bb{I}_{ B \cap \{
\bs{\tau}_2 = \iy \}}$ as $i \to \iy$.  Therefore, applying the bounded convergence theorem leads to (\ref{claim82}) and the inequality $\int_{
B \cap \{ \bs{\tau}_2 = \iy \}}  | X_{\bs{\tau}_2}  | \; d \bb{P} < \iy$.

Finally, combining (\ref{claim81}) and (\ref{claim82}) gives
\[
\lim_{i \to \iy} \int_{ B \cap \{ \bs{\tau}_2 > i \}}  X_i  \; d \bb{P} = \lim_{i \to \iy} \int_{ B \cap \{ i < \bs{\tau}_2 < \iy \}}  X_i  \; d
\bb{P} + \lim_{i \to \iy} \int_{ B \cap \{ \bs{\tau}_2 = \iy \}}  X_i  \; d \bb{P} = \int_{ B \cap \{ \bs{\tau}_2 = \iy \}}  X_{\bs{\tau}_2}  \;
d \bb{P}.
\]
This completes the proof of the lemma.

\epf

\bsk

Now we are in a position to prove the discrete-time optional stopping theorem under assumption (\ref{asmp888}).  Since $X_{\bs{\tau}_2}$ exits,
it must be true that $X_{\bs{\tau}_1}$ exists. It suffices to show that
\[
\int_{ A \cap  \{ \bs{\tau}_2 \geq \bs{\tau}_1 \}}  X_{\bs{\tau}_2} \; d \bb{P} \leq \int_{ A \cap  \{ \bs{\tau}_2 \geq \bs{\tau}_1 \}}
X_{\bs{\tau}_1} \; d \bb{P}
\]
for any $A \in \mcal{F}_{ \bs{\tau}_1 }$.  For this, in turn, it is sufficient to show that, for every $n \in \bb{Z}^+ \cup \{ \iy \}$,
\[
\int_{ A \cap  \{ \bs{\tau}_2 \geq \bs{\tau}_1 \} \cap \{ \bs{\tau}_1 = n \} }  X_{\bs{\tau}_2} \; d \bb{P} \leq \int_{ A \cap  \{ \bs{\tau}_2
\geq \bs{\tau}_1 \} \cap \{ \bs{\tau}_1 = n \} } X_{\bs{\tau}_1} \; d \bb{P}.
\]
This inequality is clearly true for $n = \iy$, because
\[
A \cap  \{ \bs{\tau}_2 \geq \bs{\tau}_1 \} \cap \{ \bs{\tau}_1 = n \} = A \cap  \{ \bs{\tau}_1 = \iy, \; \bs{\tau}_2 = \iy \} = A \cap  \{
\bs{\tau}_1 = \iy, \; \bs{\tau}_2 = \iy, \; X_{\bs{\tau}_2} = X_{\bs{\tau}_1}  \}
\]
holds for $n = \iy$.  Recall Lemma \ref{lem1a}, we have that for $i \geq n$, \be \la{inqlim}
 \int_{ B \cap \{ \bs{\tau}_2 \geq n \}}  X_n  \; d \bb{P} \geq \int_{ B \cap \{ n
\leq \bs{\tau}_2 \leq i \}} X_{\bs{\tau}_2} \; d \bb{P} + \int_{ B \cap \{ \bs{\tau}_2 > i \}}  X_i  \; d \bb{P}, \ee  with equality if $X_t$ is
a martingale.  Taking limits on the right side of (\ref{inqlim}) and making use of Lemmas \ref{lem889} and \ref{lem8899}, we have \bee \int_{ B
\cap \{ \bs{\tau}_2 \geq n \}}  X_n  \; d \bb{P} & \geq & \lim_{i \to \iy} \int_{ B \cap \{ n \leq \bs{\tau}_2 \leq i \}} X_{\bs{\tau}_2} \; d
\bb{P} + \lim_{i \to \iy} \int_{ B \cap \{ \bs{\tau}_2 > i \}}  X_i  \; d \bb{P}\\
& = &  \int_{ B \cap \{ n \leq \bs{\tau}_2 < \iy \}} X_{\bs{\tau}_2} \; d \bb{P} +  \int_{ B \cap \{ \bs{\tau}_2 = \iy \}}  X_{\bs{\tau}_2}  \;
d \bb{P} =   \int_{ B \cap \{ \bs{\tau}_2 \geq n \}} X_{\bs{\tau}_2} \; d \bb{P},  \eee  with equality if $(X_k, \mscr{F}_k)_{k \in \bb{Z}^+}$
is a martingale.  This completes the proof of the discrete-time optional stopping theorem under assumption (\ref{asmp888}).

\subsection{Proof of Discrete-Time Optional Stopping Theorem for Super-martingale with Bounded Increment}

In this section of Appendix \ref{Themmain_app},  we shall show the discrete-time optional stopping theorem for super-martingale with bounded
increment. More formally, we want to prove the following result:

{\it Let $(X_k, \mscr{F}_k)_{k \in \bb{Z}^+ }$ be a super-martingale such that there exists a constant $\vDe$ so that $| X_{k+1} - X_k | < \vDe$
almost surely for all $k \in \bb{Z}^+$.  Let $\bs{\tau}_1$ and $\bs{\tau}_2$ be two stopping times such that $\bs{\tau}_1 \leq \bs{\tau}_2$
almost surely and that there exists a constant $C$ so that $\{ \bs{\tau}_2 > k \} \subseteq \{ | X_k | < C \}$ for all $k \in \bb{Z}^+$.  Then,
$\bb{E} [ X_{ \bs{\tau}_2  } \mid \mscr{F}_{ \bs{\tau}_1 } ] \leq X_{ \bs{\tau}_1 }$ and $\bb{E} [ X_{ \bs{\tau}_2 } ] \leq \bb{E} [ X_{
\bs{\tau}_1 } ]$ almost surely, with equality if $(X_k, \mscr{F}_k)_{k \in \bb{Z}^+}$ is a martingale.}

For this purpose, it suffices to show that assumption (\ref{asmp888}) is satisfied.  This can be accomplished by proving the following lemma.

\beL

\la{lembounded} $X_{ \bs{\tau} \wedge k} $ is a UI super-martingale. Moreover, $X_{ \bs{\tau}} $ exists and $\bb{E} [ | X_{ \bs{\tau}} | ] <
\iy$.  \eeL

\bpf

For simplicity of notations, we denote $\int_F X  d \bb{P}$ by $\bb{E} [ X; F]$. Let $\Up > 0$.   Note that \bee &  &  \bb{E} [| X_{ \bs{\tau}
\wedge k} |; \; | X_{ \bs{\tau} \wedge k} | \geq \Up ] = \bb{E} [| X_{ \bs{\tau} \wedge k}  |; \; | X_{ \bs{\tau} \wedge k}  | \geq \Up, \;
\bs{\tau} > k ] + \sum_{i = 0}^{k} \bb{E} [| X_{ \bs{\tau} \wedge k} |; \; | X_{ \bs{\tau} \wedge k}
| \geq \Up, \; \bs{\tau} = i ]\\
&  &  = \bb{E} [| X_k  |; \;   | X_k | \geq \Up, \; \bs{\tau} > k ] +
 \sum_{i = 0}^{k} \bb{E} [| X_i  |; \;   | X_{ \bs{\tau}} | \geq \Up, \; \bs{\tau} = i ]\\
&  &  = \bb{E} [| X_k  |; \;   | X_k | \geq \Up, \; \bs{\tau} > k ] + \bb{E} [| X_0  |; \;   | X_{ \bs{\tau}} | \geq \Up, \; \bs{\tau} = 0 ]
 +
\sum_{i = 1}^{k} \bb{E} [| X_i  |; \;   | X_{
\bs{\tau}} | \geq \Up, \; i - 1 < \bs{\tau} = i ]\\
&  &  = \bb{E} [| X_k  |; \;   | X_k | \geq \Up, \; \bs{\tau} > k ] + \bb{E} [| X_0  |; \;   | X_0 | \geq
\Up, \; \bs{\tau} = 0 ] + \sum_{i = 1}^{k} \bb{E} [| X_i  |; \; | X_i | \geq \Up, \; i - 1 < \bs{\tau} = i ]\\
&  &  \leq \bb{E} [| X_k  |; \;  \Up \leq |X_k| < C ] + \bb{E} [| X_0  |; \;   | X_0 | \geq \Up ] + \sum_{i = 1}^{k} \bb{E} [| X_i  |; \;
  | X_i | \geq \Up, \; |X_{i-1}| < C  ]\\
&  &  \leq \bb{E} [| X_k  |; \; \Up \leq |X_k| < C ] + \bb{E} [| X_0  |; \;   | X_0 | \geq \Up ]
 + \sum_{i = 1}^{k} \bb{E} [| X_i  |; \; |X_{i-1}| +
| X_i - X_{i-1} | \geq \Up, \; |X_{i-1}| < C  ] \eee for all $k \geq 0$. By the bounded increment assumption, there exists a positive constant
$\vDe
> 0$ such that $\Pr \{  | X_i - X_{i-1} | < \vDe  \} = 1$ for $i \geq 1$.  Choose $\Up
> C + \vDe$, then $\bb{E} [| X_{ \bs{\tau} \wedge k}  |; \;   | X_{ \bs{\tau} \wedge k}  | \geq \Up ]
 \leq \bb{E} [| X_0  |; \;   | X_0 | \geq \Up ]$ for all $k \geq 0$.
 By the assumption that $\bb{E}[|X_0|] < \iy$, we have that  there exists a sufficiently large $\Up > C + \vDe$ such
that $\bb{E} [| X_0  |; \;   | X_0 | \geq \Up ] < \vep$.   Hence, $\bb{E} [ |X_{ \bs{\tau} \wedge k}|; \; |X_{ \bs{\tau} \wedge k}| \geq \Up ] <
\vep$ for all $k \geq 0$. This implies that $X_{ \bs{\tau} \wedge k} $ is a UI super-martingale.

Since $X_{ \bs{\tau} \wedge k} $ is a UI super-martingale, we have $\sup_{k} \bb{E} [ | X_{ \bs{\tau} \wedge k} | ] < \iy$. Moreover, $X_{
\bs{\tau} \wedge k} $ converges as $k \to \iy$ and $X_{ \bs{\tau}} $ exists almost surely. By Fatou's lemma, $\bb{E} [ | X_{ \bs{\tau}} | ] =
\bb{E} [ \liminf_{k \to \iy } | X_{ \bs{\tau} \wedge k} | ] \leq \liminf_{k \to \iy } \bb{E} [ | X_{ \bs{\tau} \wedge k} | ] \leq \sup_{k}
\bb{E} [ | X_{ \bs{\tau} \wedge k} | ]  < \iy$.

\epf

\subsection{Proof of Discrete-Time Optional Stopping Theorem for Non-negative Super-martingale}

In this section of Appendix \ref{Themmain_app},  we shall show the discrete-time optional stopping theorem for non-negative super-martingale.
More formally, we want to prove the following result:

{\it Let $(X_k, \mscr{F}_k)_{k \in \bb{Z}^+ }$ be a non-negative super-martingale.  Let $\bs{\tau}_1$ and $\bs{\tau}_2$ be two stopping times
such that $\bs{\tau}_1 \leq \bs{\tau}_2$ almost surely and that there exists a constant $C$ so that $\{ \bs{\tau}_2 > k \} \subseteq \{ X_k < C
\}$ for all $k \in \bb{Z}^+$.   Then, $\bb{E} [ X_{ \bs{\tau}_2  } \mid \mscr{F}_{ \bs{\tau}_1 } ] \leq X_{ \bs{\tau}_1 }$ and $\bb{E} [ X_{
\bs{\tau}_2 } ] \leq \bb{E} [ X_{ \bs{\tau}_1 } ]$ almost surely, with equality if $(X_k, \mscr{F}_k)_{k \in \bb{Z}^+}$ is a martingale.}

\bsk

 Since a non-negative super-martingale must converge, it follows that both $X_{\bs{\tau}_1}$ and $X_{\bs{\tau}_2}$ exist.  As an immediate
consequence of Lemma \ref{lem8899}, we have the following result.

\beL \la{claimtop} For any non-negative integer $n$,  $0 \leq \lim_{i \to \iy} \int_{ B \cap \{ \bs{\tau}_2 > i \}}  X_i  \; d \bb{P} = \int_{ B
\cap \{ \bs{\tau}_2 = \iy \}} X_{\bs{\tau}_2} \; d \bb{P} < \iy$.  \eeL

\bsk

We are now in a position to prove the discrete-time optional stopping theorem for non-negative super-martingale.  It suffices to show that
\[
\int_{ A \cap  \{ \bs{\tau}_2 \geq \bs{\tau}_1 \}}  X_{\bs{\tau}_2} \; d \bb{P} \leq \int_{ A \cap  \{ \bs{\tau}_2 \geq \bs{\tau}_1 \}}
X_{\bs{\tau}_1} \; d \bb{P}
\]
for any $A \in \mscr{F}_{ \bs{\tau}_1 }$. For this, in turn, it is sufficient to show that, for every $n \in \bb{Z}^+ \cup \{ \iy \}$,
\[
\int_{ A \cap  \{ \bs{\tau}_2 \geq \bs{\tau}_1 \} \cap \{ \bs{\tau}_1 = n \} }  X_{\bs{\tau}_2} \; d \bb{P} \leq \int_{ A \cap  \{ \bs{\tau}_2
\geq \bs{\tau}_1 \} \cap \{ \bs{\tau}_1 = n \} } X_{\bs{\tau}_1} \; d \bb{P}.
\]
This inequality is clearly true for $n = \iy$, because
\[
A \cap  \{ \bs{\tau}_2 \geq \bs{\tau}_1 \} \cap \{ \bs{\tau}_1 = n \} = A \cap  \{ \bs{\tau}_1 = \iy, \; \bs{\tau}_2 = \iy \} = A \cap  \{
\bs{\tau}_1 = \iy, \; \bs{\tau}_2 = \iy, \; X_{\bs{\tau}_2} = X_{\bs{\tau}_1}  \}
\]
holds for $n = \iy$.   It remains to show, for $n \in \bb{Z}^+$, \be \la{arg1}
 \int_{ B \cap  \{ \bs{\tau}_2 \geq n \}}
X_{\bs{\tau}_2} \; d \bb{P} \leq \int_{ B \cap  \{ \bs{\tau}_2 \geq n \}} X_n \; d \bb{P}.  \ee As a consequence of Lemma \ref{recall}, \be
\la{fse} \int_{ B \cap \{ n \leq \bs{\tau}_2 < \iy \}}  X_{\bs{\tau}_2} \; d \bb{P} = \bb{E} \li [ \liminf_{i \to \iy} X_{\bs{\tau}_2} \bb{I}_{
B \cap \{ n \leq \bs{\tau}_2 \leq i \}} \ri ]. \ee  Note that  \bel \bb{E} \li [  \liminf_{i \to \iy} X_{\bs{\tau}_2} \bb{I}_{ B \cap \{ n \leq
\bs{\tau}_2 \leq i \}} \ri ] & \leq & \liminf_{i \to \iy} \bb{E} [
X_{\bs{\tau}_2} \bb{I}_{ B \cap \{ n \leq \bs{\tau}_2 \leq i \}} ] \la{usefatou}\\
& = & \liminf_{i \to \iy} \int_{ B \cap \{ n \leq \bs{\tau}_2 \leq i \}}  X_{\bs{\tau}_2} \; d \bb{P} \nonumber\\
& \leq & \liminf_{i \to \iy} \li [ \int_{ B \cap \{ \bs{\tau}_2 \geq n \}} X_n  \; d \bb{P} - \int_{ B \cap \{ \bs{\tau}_2 > i
\}}  X_i  \; d \bb{P} \ri ] \la{uselem1a}\\
& = &  \int_{ B \cap \{ \bs{\tau}_2 \geq n \}} X_n  \; d \bb{P} - \limsup_{i \to \iy} \int_{ B \cap \{ \bs{\tau}_2 > i \}}  X_i  \; d \bb{P} \nonumber\\
& = &  \int_{ B \cap \{ \bs{\tau}_2 \geq n \}} X_n  \; d \bb{P} - \int_{ B \cap \{ \bs{\tau}_2 = \iy \}}  X_{\bs{\tau}_2} \; d \bb{P},
\la{useclaimtop} \eel where (\ref{usefatou}) follows from Fatou's lemma, (\ref{uselem1a}) follows from Lemma \ref{lem1a}, and
(\ref{useclaimtop}) follows from Lemma \ref{claimtop}.  By the assumption that $(X_k, \mscr{F}_k)_{k \in \bb{Z}^+}$  is a super-martingale, we
have that $X_n$ is integrable and thus  $0 \leq \int_{ B \cap \{ \bs{\tau}_2 \geq n \}} X_n \; d \bb{P} < \iy$. From Lemma \ref{claimtop}, we
know that $0 \leq \int_{ B \cap \{ \bs{\tau}_2 = \iy \}}  X_{\bs{\tau}_2} \; d \bb{P} < \iy$. It follows from (\ref{useclaimtop}) that $0 \leq
\bb{E} \li [ \liminf_{i \to \iy} X_{\bs{\tau}_2} \bb{I}_{ B \cap \{ n \leq \bs{\tau}_2 \leq i \}} \ri ] < \iy$. Combing (\ref{fse}) and
(\ref{useclaimtop}) yields
\[
0 \leq \int_{ B \cap \{ n \leq \bs{\tau}_2 < \iy \}}  X_{\bs{\tau}_2} \; d \bb{P} \leq \int_{ B \cap \{ \bs{\tau}_2 \geq n \}} X_n  \; d \bb{P}
- \int_{ B \cap \{ \bs{\tau}_2 = \iy \}}  X_{\bs{\tau}_2} \; d \bb{P} < \iy
\]
or equivalently,
\[
0 \leq \int_{ B \cap \{ \bs{\tau}_2 \geq n \}}  X_{\bs{\tau}_2} \; d \bb{P} - \int_{ B \cap \{ \bs{\tau}_2 = \iy \}}  X_{\bs{\tau}_2} \; d
\bb{P} \leq \int_{ B \cap \{ \bs{\tau}_2 \geq n \}} X_n  \; d \bb{P} - \int_{ B \cap \{ \bs{\tau}_2 = \iy \}}  X_{\bs{\tau}_2} \; d \bb{P} <
\iy,
\]
which implies (\ref{arg1}).  So, we have established that $\bb{E} [ X_{\bs{\tau}_2}  \mid  \mscr{F}_{\bs{\tau}_1}  ] \leq X_{\bs{\tau}_1}$ and
$\bb{E} [ X_{\bs{\tau}_2} ] \leq \bb{E} [ X_{\bs{\tau}_1} ]$ for the case that $(X_k, \mscr{F}_k)_{k \in \bb{Z}^+}$ is a super-martingale.  It
follows that $0 \leq \bb{E} [ X_{\bs{\tau}_2} ] \leq \bb{E} [ X_0 ] < \iy$.   Applying Theorem \ref{Themmain}, we have  $\bb{E} [
X_{\bs{\tau}_2}  \mid \mscr{F}_{\bs{\tau}_1}  ] = X_{\bs{\tau}_1}$ and $\bb{E} [ X_{\bs{\tau}_2} ] = \bb{E} [ X_{\bs{\tau}_1} ]$ for the case
that $(X_k, \mscr{F}_k)_{k \in \bb{Z}^+}$  is a martingale.  Thus, we have established the discrete-time optional stopping theorem for
non-negative super-martingale.

So, we have completed the proof of Theorem \ref{Themmain}.

\section{Proof of Theorem \ref{Them3}}  \la{Them3_app}

Recall the assumption that there exist  constants $\de$ and $\vDe$ such that $| X_{t^\prime} - X_t | < \vDe$ almost surely provided that $|
t^\prime - t | < \de$.  For such $\de$, define $\nu = \lc \log_2 \f{1}{\de} \rc$ and $\bb{T}_n = \{ k 2^{-(n+ \nu)}: k \in \bb{Z}^+ \}$ for $n
\in \bb{N}$.  Define \bee &  & \bs{\ro}_n = \inf \{ t \in \bb{T}_n: X_t \notin \mscr{D}_t \}, \qqu \; \bs{\vro}_n = \inf \{ t \in  \bb{T}_n: X_t
\notin \mcal{D}_t \}, \\
&  & X_{\bs{\ro}_n} = \lim_{t \to \iy \atop{ t \in \bb{T}_n } } X_{ \bs{\ro}_n \wedge t}, \qqu \qqu \qqu \; X_{\bs{\vro}_n} = \lim_{t \to \iy
\atop{ t \in \bb{T}_n } } X_{ \bs{\vro}_n \wedge t},\\
&  & S_n = 2^{n + \nu} \bs{\ro}_n, \qqu \qqu \qqu \qqu \; \; T_n = 2^{n + \nu} \bs{\vro}_n \eee for $n \in \bb{N}$.  Define \[ t_{k, n} = k
2^{-(n+ \nu)}, \qqu Y_k^n = X_{t_{k, n}}, \qqu \mcal{F}_k^n = \mscr{F}_{t_{k, n}} \] for $k \in \bb{Z}^+$ and $n \in \bb{N}$. Then, the
following statements are true.

\noindent (a): For $n \in \bb{N}$, $(Y_k^n, \mcal{F}_k^n)_{k \in \bb{Z}^+}$ is a discrete-time super-martingale.

\noindent (b): For $k \in \bb{Z}^+$ and $n \in \bb{N}$, the inequality $t_{k + 1, n} - t_{k, n} < \de$ holds and consequently, $| Y_{k+1}^n -
Y_k^n | < \vDe$ almost surely.

\noindent (c):  $S_n \leq T_n$ and $\{ T_n > k \} \subseteq \{ |Y_k^n | > C \}$ for $n \in \bb{N}$.

\bsk

Define \[ Y_{S_n}^n = \lim_{k \to \iy} Y_{S_n \wedge k}^n, \qqu Y_{T_n}^n = \lim_{k \to \iy} Y_{T_n \wedge k}^n
\]
for $n \in \bb{N}$.  Note that $S_n$ and $T_n$ are stopping times non-increasing with respect to $n \in \bb{N}$.   To complete the proof of the
theorem, we need the following results.

\beL \la{lem999}

(I) $\bs{\ro}_n$ and $\bs{\vro}_n$ are stopping times non-increasing with respect to $n \in \bb{N}$.

(II) $\bs{\ro}_n \geq \bs{\tau}_1, \; \bs{\vro}_n \geq \bs{\tau}_2$ and $\bs{\ro}_n \leq \bs{\vro}_n$ for all $n \in \bb{N}$.

(III) $\lim_{n \to \iy} \bs{\ro}_n = \bs{\tau}_1$ and $\lim_{n \to \iy} \bs{\vro}_n = \bs{\tau}_2$.

(IV)  For all $n \in \bb{N}$,  $X_{\bs{\ro}_n}$ and $X_{\bs{\vro}_n}$ exist almost surely.

(V) $X_{\bs{\tau}_1}$ and $X_{\bs{\tau}_2}$  exist almost surely.

(VI) As $n$ tends to infinity, $X_{\bs{\ro}_n}$ and $X_{\bs{\vro}_n}$ converge to $X_{\bs{\tau}_1}$ and $X_{\bs{\tau}_2}$ respectively and
almost surely.

\eeL

\bpf

Statements (I) -- (III) are obviously true.  We shall show statements (IV), (V) and (VI).

\bed

\item [{\it Proof of Statement (IV)}]:  Consider the existence of $X_{\bs{\ro}_n}$ for $n \in \bb{N}$.   From Lemma \ref{lembounded} and statements
(a), (b) and (c) appeared before Lemma \ref{lem999}, we know that for every $n \in \bb{N}$, $(Y_{S_n \wedge k}^n, \mcal{F}_k^n)_{k \in
\bb{Z}^+}$ is a discrete-time UI martingale and it follows that, almost surely, $\lim_{k \to \iy} Y_{S_n \wedge k}^n$ exists and is finite.  By
the definitions of $\bb{T}_n, \; \bs{\ro}_n, \; S_n$ and $\{Y_k^n\}$, we have that $\lim_{t \to \iy \atop{ t \in \bb{T}_n} } X_{ \bs{\ro}_n
\wedge t} = \lim_{k \to \iy} Y_{S_n \wedge k}^n$ almost surely for all $n \in \bb{N}$.  Since $X_{ \bs{\ro}_n }$ is defined as $\lim_{t \to \iy
\atop{ t \in \bb{T}_n} } X_{ \bs{\ro}_n \wedge t}$,  it follows that $X_{\bs{\ro}_n}$ exists almost surely for all $n \in \bb{N}$. In a similar
manner, the existence of $X_{\bs{\vro}_n}$ can be established for $n \in \bb{N}$.

\item [{\it Proof of Statement (V)}]: Consider the existence of $X_{\bs{\tau}_1}$.
Let $n \in \bb{N}$ be fixed and let $\om \in \Om$ with $\bs{\tau}_1 (\om) = \iy$.  From the proof of Statement (IV), we know that  the limit
$\lim_{t \to \iy \atop{ t \in \bb{T}_n} } X_{ \bs{\ro}_n \wedge t}$ exists.  Since $\bs{\ro}_n  \geq \bs{\tau}_1 $,  we have $\bs{\ro}_n (\om) =
\iy$.  This implies that $\lim_{t \to \iy \atop{ t \in \bb{T}_n} } X_{ t} (\om)$ exists.   We claim that the limit $\lim_{t \to \iy} X_t (\om)$
exists and is equal to $\lim_{t \to \iy \atop{ t \in \bb{T}_n} } X_{ t} (\om)$.  Suppose, to get a contradiction, that $\lim_{t \to \iy  } X_{t}
(\om)$ does not exist. Then,  there exist an $\vep > 0$ and a sequence $\{ t_i \}_{i = 1}^\iy$ with $t_1 < t_2 < t_3 < \cd$ and $t_i \to \iy$ as
$i \to \iy$, such that $| X_{t_i} - c | > \vep$ for $i \geq 1$, where $c$ denotes $\lim_{t \to \iy \atop{ t \in \bb{T}_n} } X_{ t} (\om)$.
Define $\bb{W}_n = \bb{T}_n \cup \{t_i\}_{i = 1}^\iy$. That is,  the sequence $\{t_i\}_{i=1}^\iy$ is added to $\bb{T}_n$ to form a new sequence
$\bb{W}_n$. Define $\bs{\mu}_n = \inf \{ t \in \bb{W}_n: X_t \notin \mscr{D}_t \}$.  By the same argument as that for proving the existence of
$\lim_{t \to \iy \atop{ t \in \bb{T}_n} } X_{ \bs{\ro}_n \wedge t}$, we can show that  $\lim_{t \to \iy \atop{ t \in \bb{W}_n} } X_{ \bs{\mu}_n
\wedge t}$ exists almost surely.  Observing that $\bs{\mu}_n  \geq \bs{\tau}_1 $,  we have $\bs{\mu}_n (\om) = \iy$. Therefore, $\lim_{t \to \iy
\atop{ t \in \bb{W}_n} } X_{ t} (\om)$ exists.  Since $\bb{T}_n \subseteq \bb{W}_n$, it must be true that $\lim_{t \to \iy \atop{ t \in
\bb{W}_n} } X_{ t} (\om) = c$.  Since $\{t_i\}_{i=1}^\iy  \subseteq \bb{W}_n$, it follows that $\lim_{i \to \iy} X_{t_i} (\om)$ exists and is
equal to $c$. This contradicts to the assumption that $| X_{t_i} - c | > \vep$ for $i \geq 1$. Thus, we have established the claim that $\lim_{t
\to \iy} X_t (\om)$ exists and is equal to $\lim_{t \to \iy \atop{ t \in \bb{T}_n} } X_{ t} (\om)$ for $\om \in \Om$ with $\bs{\tau}_1 (\om) =
\iy$.  Since $X_{ \bs{\tau}_1 }$ is defined as $\lim_{t \to \iy } X_{ \bs{\tau}_1 \wedge t}$,  it follows that $X_{\bs{\tau}_1}$ exists almost
surely. In a similar manner, the existence of $X_{\bs{\tau}_2}$ can be established.

\item [{\it Proof of Statement (IV)}]: Consider the convergence of $( X_{\bs{\ro}_n} )_{n \in \bb{N}}$.  Recall the established fact that
$X_{ \bs{\ro}_n } = \lim_{t \to \iy \atop{ t \in \bb{T}_n} } X_{ \bs{\ro}_n \wedge t}$ exists almost surely for all $n \in \bb{N}$.  Let $\om
\in \Om$ with $\bs{\tau}_1 (\om) = \iy$.  Since $ \bs{\ro}_n \geq \bs{\tau}_1$, we have $\bs{\ro}_n (\om) = \iy$ for all $n \in \bb{N}$.  It
follows that $X_{ \bs{\ro}_n } (\om) = \lim_{t \to \iy \atop{ t \in \bb{T}_n} } X_{t} (\om) = \lim_{t \to \iy } X_{t} (\om)$ for all $n \in
\bb{N}$. Therefore, \be \la{comba899}
 \lim_{n \to \iy} X_{ \bs{\ro}_n } (\om) = \lim_{t \to \iy } X_{t} (\om) \qu
\tx{for $\om \in \Om$ with $\bs{\tau}_1 (\om) = \iy$}. \ee Now let $\om \in \Om$ with $\bs{\tau}_1 (\om) < \iy$.  Since $\lim_{n \to \iy}
\bs{\ro}_n = \bs{\tau}_1$, we have that $\bs{\ro}_n (\om) < \iy$ for sufficiently large $n \in \bb{N}$.  This implies that $X_{ \bs{\ro}_n }
(\om) = X_{ \bs{\ro}_n (\om) } (\om)$ for sufficiently large $n \in \bb{N}$.  Since $\lim_{n \to \iy} \bs{\ro}_n (\om) = \bs{\tau}_1 (\om)$ and
$( X_t )_{t \in \bb{R}^+}$ is a right-continuous process, we have that \be \la{combb899} \lim_{n \to \iy} X_{ \bs{\ro}_n } (\om) =
X_{\bs{\tau}_1 (\om)} (\om) \qu \tx{for $\om \in \Om$ with $\bs{\tau}_1 (\om) < \iy$}. \ee Making use of (\ref{comba899}), (\ref{combb899}) and
the definition of $X_{\bs{\tau}_1}$, we have that $X_{ \bs{\ro}_n }$ converges to $X_{\bs{\tau}_1}$ almost surely.  Similarly, we can show that
$X_{ \bs{\vro}_n }$ converges to $X_{\bs{\tau}_2}$ almost surely.

\eed

\epf

\beL \la{lem998}

As $n$ tends to infinity, $X_{\bs{\ro}_n}$ and $X_{\bs{\vro}_n}$ converge to $X_{\bs{\tau}_1}$ and $X_{\bs{\tau}_2}$ respectively in $L^1$.

\eeL

\bpf

Let $\{ \mscr{G}_{-n}:  n \in  \bb{N} \cup \{ \iy \} \}$ be a collection of sub-$\si$-algebras of $\mscr{F}$ with $\mscr{G}_{-n} =
\mscr{F}_{\bs{\ro}_n}$ for $n \in \bb{N}$ and $\mscr{G}_{- \iy} = \bigcap_{k\in \bb{N}} \mscr{G}_{- k}$.  Define $Z_{-n} = X_{ \bs{\ro}_n }$ for
$n \in \bb{N}$.  Then, $Z_{-n} = Y_{S_n}^n$ for $n \in \bb{N}$. According to Theorem \ref{Themmain}, we have that $\bb{E} [ Y_{S_n}^n  \mid
\mcal{F}_{S_{n+1}}^{n+1} ] \leq Y_{S_{n+1}}^{n+1}$ almost surely for $n \in \bb{N}$.  Since $Y_{S_n}^n = X_{  \bs{\ro}_n }$ and
$\mcal{F}_{S_n}^n = \mscr{F}_{\bs{\ro}_n}$ for $n \in \bb{N}$, we have that $\bb{E} [ X_{ \bs{\ro}_n } \mid \mscr{F}_{\bs{\ro}_{n+1}} ] \leq
X_{\bs{\ro}_{n+1}}$ almost surely for $n \in \bb{N}$.  This implies that $\bb{E} [ Z_{-n} \mid \mscr{G}_{-(n+1)} ] \leq Z_{-(n+1)}$ almost
surely.  Hence, $\{ Z_{-n}, \; n \in  \bb{N} \}$ is a supermartingale relative to $\{ \mscr{G}_{-n}:  n \in \bb{N} \cup \{ \iy \} \}$ in the
context of L$\acute{e}$vy-Doob Downward Theorem (see, e.g., \cite[page 148--149]{Levy-Doob}).  Moreover, from Theorem \ref{Themmain}, we have
that $\bb{E} [ Y_{S_n}^n ] \leq \bb{E} [ Y_{0, n} ]$ almost surely for $n \in \bb{N}$.  Since $Y_{S_n}^n = X_{  \bs{\ro}_n }$ and $Y_{0, n} =
X_0$, we have $\bb{E} [ X_{ \bs{\ro}_n } ] \leq \bb{E} [ X_0] < \iy$, which implies that $\sup_{n \in \bb{N}} \bb{E} [ Z_{-n} ] = \sup_{n \in
\bb{N}} \bb{E} [ X_{ \bs{\ro}_n } ] < \iy$. Therefore, it follows from L$\acute{e}$vy-Doob Downward Theorem that $\{ Z_{-n}, \; n \in \bb{N} \}$
is uniformly integrable and that the limit $Z_{-\iy} \DEF \lim_{n \to \iy} Z_{-n}$ exists almost surely and $\bb{E}[ | Z_{-\iy} | ] < \iy$. From
Lemma \ref{lem999}, we know that $Z_{-\iy} = X_{ \bs{\tau}_1}$.  Since $\{ Z_{-n}, \; n \in \bb{N} \}$ is uniformly integrable and converges to
$X_{\bs{\tau}_1}$ almost surely, it follows that $X_{ \bs{\ro}_n } = Z_{-n} \to X_{ \bs{\tau}_1 }$ in $L^1$. Similarly, we have that $X_{
\bs{\vro}_n } \to X_{ \bs{\tau}_2 }$ in $L^1$.

\epf

\bsk

We are now in a position to prove the theorem.  We follow the classical argument for extending the optional stopping theorem from discrete UI
martingale to continuous-time UI martingale.  According to Theorem \ref{Themmain}, we have $\bb{E} [ Y_{T_n}^n \mid \mcal{F}_{S_n}^n ] \leq
Y_{S_n}^n$ almost surely for $n \in \bb{N}$. Since $Y_{S_n}^n = X_{  \bs{\ro}_n }, \; Y_{T_n}^n = X_{  \bs{\vro}_n }$ and $\mcal{F}_{S_n}^n =
\mscr{F}_{\bs{\ro}_n}$, we have \be \bb{E} [ X_{  \bs{\vro}_n } \mid \mscr{F}_{\bs{\ro}_n} ] \leq X_{\bs{\ro}_n} \la{qq18} \ee almost surely for
$n \in \bb{N}$. Since $\bs{\tau}_1 \leq \bs{\ro}_n$, it holds that $\mscr{F}_{\bs{\tau}_1} \subseteq \mscr{F}_{\bs{\ro}_n}$. It follows from
(\ref{qq18}) and the tower property that \be \la{qq28}
 \bb{E} [ X_{  \bs{\vro}_n } \mid \mscr{F}_{\bs{\tau}_1} ] =  \bb{E} [ \bb{E} [ X_{  \bs{\vro}_n } \mid \mscr{F}_{\bs{\ro}_n} ] \mid \mscr{F}_{\bs{\tau}_1}
 ] \leq \bb{E} [ X_{  \bs{\ro}_n } \mid \mscr{F}_{\bs{\tau}_1} ]
\ee almost surely.  Let $\mscr{E} \in \mscr{F}_{ \bs{\tau}_1 }$. It follows from (\ref{qq28}) that $\int_\mscr{E} X_{ \bs{\vro}_n } d \bb{P}
\leq \int_\mscr{E} X_{ \bs{\ro}_n } d \bb{P}$.  Invoking Lemma \ref{lem998}, we have $X_{ \bs{\ro}_n } \to X_{ \bs{\tau}_1 }$ in $L^1$ and
consequently,
\[
\li | \int_\mscr{E} X_{  \bs{\ro}_n } d \bb{P} - \int_\mscr{E} X_{  \bs{\tau}_1 } d \bb{P}  \ri |  \leq \int_\mscr{E} | X_{  \bs{\ro}_n } - X_{
\bs{\tau}_1 } | \; d \bb{P} \leq \bb{E} [ | X_{  \bs{\ro}_n } - X_{ \bs{\tau}_1 } | ] \to 0
\]
as $n \to 0$.  This implies that $ \int_\mscr{E} X_{  \bs{\ro}_n } d \bb{P} \to  \int_\mscr{E} X_{ \bs{\tau}_1 } d \bb{P}$ as $n \to 0$.
Similarly,  $ \int_\mscr{E} X_{ \bs{\vro}_n } d \bb{P} \to  \int_\mscr{E} X_{  \bs{\tau}_2 } d \bb{P}$ as $n \to 0$.  Therefore,  $\int_\mscr{E}
X_{ \bs{\tau}_2 } d \bb{P} = \lim_{n \to \iy} \int_\mscr{E} X_{ \bs{\vro}_n } d \bb{P} \leq \lim_{n \to \iy} \int_\mscr{E} X_{ \bs{\ro}_n } d
\bb{P} = \int_\mscr{E} X_{ \bs{\tau}_1 } d \bb{P}$. Since the argument holds for arbitrary $\mscr{E} \in \mscr{F}_{ \bs{\tau}_1 }$, the proof of
the theorem is thus completed.

\section{Proof of Theorem \ref{Them1}}  \la{Them1_app}

Since $(X_t, \mscr{F}_t )_{t \in \bb{R}^+}$ is a nonnegative, right-continuous supermartingale, the limit $\lim_{t \to \iy} X_t$ exists almost
surely. As a consequence of this fact and  the definition  that $X_{\bs{\tau}_i} = \lim_{t \to \iy} X_{\bs{\tau}_i \wedge t}, \; i = 1, 2$, it
must be true that both $X_{\bs{\tau}_1}$ and $X_{\bs{\tau}_2}$ exist almost surely.  For $n \in \bb{N}$, define \bee &  & \bs{\ro}_n = \inf \{ t
\in \bb{S}_n: X_t \notin \mscr{D}_t \}, \qqu \; \bs{\vro}_n = \inf \{ t \in  \bb{S}_n: X_t
\notin \mcal{D}_t \}, \\
&  & X_{\bs{\ro}_n} = \lim_{t \to \iy \atop{ t \in \bb{S}_n } } X_{ \bs{\ro}_n \wedge t}, \qqu \qqu \qqu X_{\bs{\vro}_n} = \lim_{t \to \iy
\atop{ t\in \bb{S}_n } } X_{ \bs{\vro}_n \wedge t},\\
&  & S_n = 2^{n + \nu} \bs{\ro}_n, \qqu \qqu \qqu \qqu \; \; T_n = 2^{n + \nu} \bs{\vro}_n. \eee  By the same argument as that for the existence
of $X_{\bs{\tau}_1}$ and $X_{\bs{\tau}_2}$, we have that $X_{\bs{\ro}_n}$ and $X_{\bs{\ro}_n}$ exist almost surely for all $n \in \bb{N}$.
Define
\[ t_{k, n} = k 2^{-n}, \qqu Y_k^n = X_{t_{k, n}}, \qqu \mcal{F}_k^n = \mscr{F}_{t_{k, n}} \] for $k \in \bb{Z}^+$ and $n \in \bb{N}$. Then, the
following statements are true:

\noindent (a): For $n \in \bb{N}$, $(Y_k^n, \mcal{F}_k^n)_{k \in \bb{Z}^+}$ is a discrete-time non-negative super-martingale.

\noindent (b):  $S_n \leq T_n$ and $\{ T_n > k \} \subseteq \{ |Y_k^n | > C \}$ for $n \in \bb{N}$.

\bsk

Define \[ Y_{S_n}^n = \lim_{k \to \iy} Y_{S_n \wedge k}^n, \qqu Y_{T_n}^n = \lim_{k \to \iy} Y_{T_n \wedge k}^n
\]
for $n \in \bb{N}$.  Note that $S_n$ and $T_n$ are stopping times non-increasing with respect to $n \in \bb{N}$.  To complete the proof of the
theorem, we need to use the following results.

\beL

\la{lem666a}

(I) $\bs{\ro}_n$ and $\bs{\vro}_n$ are stopping times non-increasing with respect to $n \in \bb{N}$.

(II) $\bs{\ro}_n \geq \bs{\tau}_1, \; \bs{\vro}_n \geq \bs{\tau}_2$ and $\bs{\ro}_n \leq \bs{\vro}_n$ for all $n \in \bb{N}$.

(III) $\lim_{n \to \iy} \bs{\ro}_n = \bs{\tau}_1$ and $\lim_{n \to \iy} \bs{\vro}_n = \bs{\tau}_2$.

(IV) As $n$ tends to infinity, $X_{\bs{\ro}_n}$ and $X_{\bs{\vro}_n}$ converge to $X_{\bs{\tau}_1}$ and $X_{\bs{\tau}_2}$ respectively and
almost surely.

\eeL

\bpf

Statements (I) -- (III) are obvious from the definition.  We shall show statement (IV).  Consider the convergence of $( X_{\bs{\ro}_n} )_{n \in
\bb{N}}$. Since $(X_t, \; \mscr{F}_t)_{t \in \bb{R}^+}$ is a non-negative supermartingale, the limit $\lim_{t \to \iy} X_{t} (\om)$ must exist
for $\om \in \Om$.  Let $\om \in \Om$ with $\bs{\tau}_1 (\om) = \iy$.  Since $ \bs{\ro}_n \geq \bs{\tau}_1$, we have $\bs{\ro}_n (\om) = \iy$
for all $n \in \bb{N}$.  It follows that $X_{ \bs{\ro}_n } (\om) = \lim_{t \to \iy \atop{ t \in \bb{S}_n} } X_{t} (\om) = \lim_{t \to \iy }
X_{t} (\om)$ for all $n \in \bb{N}$. Therefore, \be \la{comba8}
 \lim_{n \to \iy} X_{ \bs{\ro}_n } (\om) = \lim_{t \to \iy } X_{t} (\om) \qu
\tx{for $\om \in \Om$ with $\bs{\tau}_1 (\om) = \iy$}. \ee Now let $\om \in \Om$ with $\bs{\tau}_1 (\om) < \iy$.  Since $\lim_{n \to \iy}
\bs{\ro}_n = \bs{\tau}_1$, we have that $\bs{\ro}_n (\om) < \iy$ for sufficiently large $n \in \bb{N}$.  This implies that $X_{ \bs{\ro}_n }
(\om) = X_{ \bs{\ro}_n (\om) } (\om)$ for sufficiently large $n \in \bb{N}$.  Since $\lim_{n \to \iy} \bs{\ro}_n (\om) = \bs{\tau}_1 (\om)$ and
$( X_t )_{t \in \bb{R}^+}$ is a right-continuous process, we have that \be \la{combb8} \lim_{n \to \iy} X_{ \bs{\ro}_n } (\om) =  X_{\bs{\tau}_1
(\om)} (\om) \qu \tx{for $\om \in \Om$ with $\bs{\tau}_1 (\om) < \iy$}. \ee Making use of (\ref{comba8}), (\ref{combb8}) and the definition of
$X_{\bs{\tau}_1}$, we have that $X_{ \bs{\ro}_n }$ converges to $X_{\bs{\tau}_1}$ almost surely.  Similarly, we can show that $X_{ \bs{\vro}_n
}$ converges to $X_{\bs{\tau}_2}$ almost surely.

\epf

Making use of Theorem \ref{Themmain}, Lemma \ref{lem666a} and a similar technique as that for proving Lemma \ref{lem998}, we can show the
following result.

\beL

\la{lem868}

As $n$ tends to infinity, $X_{\bs{\ro}_n}$ and $X_{\bs{\vro}_n}$ converge to $X_{\bs{\tau}_1}$ and $X_{\bs{\tau}_2}$ respectively in $L^1$.

\eeL

Finally, the proof of Theorem \ref{Them1} can be completed by making use of Theorem \ref{Themmain}, Lemma \ref{lem868} and a similar technique
as that for proving Theorem \ref{Them3}.

\section{Proof of Theorem \ref{Them6UI}}  \la{Them6UI_app}

First, let $\ga > c$ and consider $\Pr \li \{  \sup_{t \geq 0} X_t \geq \ga  \ri \}$.  Define $\bs{\tau} = \inf \{  t \in [0, \iy): X_t \geq \ga
\}$.  Then, $\bs{\tau}$ is a stopping time.  Define $X_{\bs{\tau}}$ such that for $\om \in \Om$,
\[
X_{\bs{\tau}} (\om) =  \bec X_{\bs{\tau} (\om)} (\om)  & \tx{if} \; \bs{\tau} (\om) < \iy, \\
c  & \tx{if} \; \bs{\tau} (\om) = \iy  \eec
\]
Since $(X_t, \mscr{F}_t)_{t \in \bb{R}^+}$ is a right-continuous UI martingale which converges to $c$, we have $\bb{E} | X_{\bs{\tau}} | < \iy$
and $\bb{E} [ X_{\bs{\tau}} ] = \bb{E} [ X_0 ] = c$. It follows that {\small \bee \ga \Pr \{ \bs{\tau} < \iy \}  =  \int_{ \{ \bs{\tau} < \iy \}
} \ga d \bb{P} \leq \int_{ \{ \bs{\tau} < \iy \} } X_{\bs{\tau}}  d \bb{P}  =  \bb{E} [ X_{\bs{\tau}} ]  - \bb{E} \li [ X_{\bs{\tau}} \;
\bb{I}_{\{ \bs{\tau} = \iy \}} \ri ] =  c  - \bb{E} \li [ c \; \bb{I}_{\{ \bs{\tau} = \iy \}} \ri ]
  =  c \Pr \{ \bs{\tau} < \iy \}, \eee}  which implies that $\Pr \{
\bs{\tau} < \iy \} = 0$.  Since $X_t \to c < \ga$ almost surely, we have $\Pr \li \{ \limsup_{t \geq 0} X_t \geq \ga  \ri \} = 0$, which implies
that $\Pr \li \{  \sup_{t \geq 0} X_t \geq \ga  \ri \} = \Pr \{ \bs{\tau} < \iy \}$.  Therefore, \be \la{con89a}
 \Pr \li \{  \sup_{t \geq 0} X_t \geq \ga  \ri \} = 0 \qu \tx{for $\ga > c$}.
\ee Now let $\ga < c$ and consider $\Pr \li \{  \inf_{t \geq 0} X_t \leq \ga  \ri \}$. Making use of (\ref{con89a}) and the observation that $(
- X_t, \mscr{F}_t)_{t \in \bb{R}^+}$is a right-continuous UI martingale which converges to $-c$, we have that $\Pr \li \{ \sup_{t \geq 0} (-
X_t) \geq (- \ga )  \ri \} = 0$ for $- \ga > - c$, which implies that  \be \la{con89b}
 \Pr \li \{  \inf_{t \geq 0} X_t \leq \ga  \ri \} = 0 \qu \tx{for $\ga < c$}.
\ee Combining (\ref{con89a}) and (\ref{con89b}) gives $\Pr \li \{  X_t = c \; \tx{for all} \; t \in [0, \iy)  \ri \} = 1$.  This completes the
proof of the theorem.

\section{Proof of Theorem \ref{Them5}}  \la{Them5_app}

Let $\ga > c$ and consider $\Pr \li \{  \sup_{t \geq 0} X_t \geq \ga  \ri \}$.  Define stopping time $\bs{\tau} = \inf \{  t \in [0, \iy): X_t
\geq \ga  \}$.  By definition of $X_{\bs{\tau}}$,
\[
X_{\bs{\tau}} (\om) =  \bec X_{\bs{\tau} (\om)} (\om)  & \tx{if} \; \bs{\tau} (\om) < \iy \\
c  & \tx{if} \; \bs{\tau} (\om) = \iy  \eec
\]
for $\om \in \Om$.  For simplicity of notations, let $\mu = \bb{E} [ X_0 ]$.  Since $X_t \to c < \ga$ almost surely, we have $\Pr \li \{
\limsup_{t \geq 0} X_t \geq \ga  \ri \} = 0$ and thus $\Pr \li \{  \sup_{t \geq 0} X_t \geq \ga  \ri \} = \Pr \{ \bs{\tau} < \iy \}$.  Since
$(X_t, \mscr{F}_t)_{t \in \bb{R}^+}$ is a right-continuous, non-negative super-martingale, it follows from Theorem \ref{Them1} that $\bb{E} [
X_{\bs{\tau}} ] \leq \bb{E} [ X_0 ] = \mu$. Therefore,  \bee &  & \ga \Pr \{ \bs{\tau} < \iy \} =  \int_{ \{ \bs{\tau} < \iy \} } \ga d \bb{P}
\leq \int_{ \{ \bs{\tau} < \iy \} } X_{\bs{\tau}} d \bb{P}
= \bb{E} [ X_{\bs{\tau}} ]  - \bb{E} \li [ X_{\bs{\tau}} \; \bb{I}_{\{ \bs{\tau} = \iy \}} \ri ]\\
&  & \leq  \mu  - \bb{E} \li [ X_{\bs{\tau}} \; \bb{I}_{\{ \bs{\tau} = \iy \}} \ri ] = \mu  - \bb{E} \li [ c \; \bb{I}_{\{ \bs{\tau} = \iy \}}
\ri ] = \mu - c (1 - \Pr \{ \bs{\tau} < \iy \} ). \eee  So, we have established the inequality $\ga \Pr \{ \bs{\tau} < \iy \} \leq \mu - c \li (
1 - \Pr \{ \bs{\tau} < \iy \} \ri )$.  Since $\ga > c$, solving this inequality with respect to $\Pr \{ \bs{\tau} < \iy \}$ yields {\small $\Pr
\{ \bs{\tau} < \iy \} \leq \f{ \mu - c }{\ga - c}$}.  It follows that  {\small $\Pr \li \{ \sup_{t \geq 0} X_t \geq \ga \ri \} \leq \f{ \mu - c
}{\ga - c}$} for $\ga > c$,  which implies that $c \leq \mu$.

Now consider $\Pr \li \{  \sup_{t \geq 0} X_t \geq \ga  \ri \}$ under additional assumption that $(X_t, \mscr{F}_t)_{t \in \bb{R}^+}$ is a
continuous martingale. Since $(X_t, \mscr{F}_t)_{t \in \bb{R}^+}$ is a continuous martingale, we have $\ga \Pr \{ \bs{\tau} < \iy \} =  \int_{
\{ \bs{\tau} < \iy \} } X_{\bs{\tau}} d \bb{P} = \bb{E} [ X_{\bs{\tau}} ]  - \bb{E} \li [ X_{\bs{\tau}} \; \bb{I}_{\{ \bs{\tau} = \iy \}} \ri ]
=  \mu  - \bb{E} \li [ X_{\bs{\tau}} \; \bb{I}_{\{ \bs{\tau} = \iy \}} \ri ] = \mu - c (1 - \Pr \{ \bs{\tau} < \iy \} )$ for $\ga > c$.
Consequently, {\small $\Pr \li \{ \sup_{t \geq 0} X_t \geq \ga  \ri \} = \Pr \{ \bs{\tau} < \iy \} = \f{ \mu - c }{\ga - c}$} for $\ga > c$.
Recalling that $c \leq \mu = \bb{E} [ X_0]$,  we have that {\small $\Pr \li \{ \sup_{t \geq 0} X_t \geq \ga  \ri \} = \f{ \bb{E} [ X_0] - c
}{\ga - c}$} for $\ga > \bb{E} [ X_0]$.  It remains to show that {\small $\Pr \li \{ \sup_{t \geq 0} X_t \geq \ga  \ri \} = 1$} for $\ga \leq
\bb{E} [ X_0]$.

 We claim that $\Pr \li \{ \sup_{t \geq 0} X_t \geq \mu \ri \} = 1$. In the case of $c = \mu$, if $\Pr \li \{ \sup_{t \geq 0} X_t \geq \mu
\ri \} < 1$, then there exists $\ep > 0$ such that $\Pr \{ X_t < \mu - \ep \; \tx{for all} \; t \geq 0 \} > 0$, which contradicts to the fact
that $X_t \to \mu = c$ almost surely.  In the cases of $c < \mu$, by the established result, $\Pr \li \{  \sup_{t \geq 0} X_t \geq \mu \ri \} =
\f{ \mu - c }{\mu - c} = 1$. This proves the claim.  Consequently,   we have $1 \geq \Pr \li \{  \sup_{t \geq 0} X_t \geq \ga  \ri \} \geq \Pr
\li \{ \sup_{t \geq 0} X_t \geq \mu  \ri \} = 1$ for $\ga \leq \mu = \bb{E} [ X_0]$.   This completes the proof of the theorem.

\section{Proof of Corollary \ref{expexact}}  \la{expexact_app}

By the assumption of the theorem, it can be readily shown that $(Y_t, \mscr{F}_t)_{t \in \bb{R}^+}$ is a non-negative martingale. It follows
that $(Y_t)_{t \in \bb{R}^+}$ converges almost surely. We claim that for $s \in (0, b)$, $\lim_{t \to \iy} Y_t = 0$ almost surely.   To show
this claim, it suffices to show that for $s \in (0, b)$, $\lim_{n \to \iy } Y_n = 0$ almost surely, where the limit is taken under the
constraint that $n \in \bb{N}$.  Let $\ga
> 0$ and $s \in (0, b)$. For $n \in \bb{N}$ and $\se \in (0, s)$,  we have $\bb{E} [ \exp ( \se X_n ) ] \leq \exp ( \varphi(\se) \mcal{V}_n  )$
and by Makov inequality,
\bee &  & \Pr \{  Y_n \geq \ga \} = \Pr \li \{  \exp ( \se X_n ) \geq \exp \li ( \f{\se}{s} [  \ln \ga  + \varphi(s) \mcal{V}_n  ] \ri ) \ri  \}\\
&   & \leq \f{ \exp ( \varphi(\se) \mcal{V}_n  ) }{  \exp \li ( \f{\se}{s} [  \ln \ga  + \varphi(s) \mcal{V}_n  ] \ri )  } = \li ( \f{1}{\ga}
\ri )^{ \se \sh s} \exp \li (  \li [ \f{ \varphi(\se) }{ \se } -  \f{ \varphi(s) }{ s } \ri ] \se \mcal{V}_n \ri ). \eee   By the assumption
that $\liminf_{n \to \iy}  (\mcal{V}_{n+1} - \mcal{V}_n ) > 0$, we have that $\mcal{V}_n > 0$ for large enough $n \in \bb{N}$.  Since
$\varphi(s) \mcal{V}_n$ is a convex function of $s$, it follows that $\varphi(s)$ is a convex function, which implies that $\f{ \varphi(\se) }{
\se } -  \f{ \varphi(s) }{ s } < 0$.  Since $\liminf_{n \to \iy}  (\mcal{V}_{n+1} - \mcal{V}_n ) > 0$, there exists a real number $d > 0$ and an
$m \in \bb{N}$ such that $\mcal{V}_{n+1} - \mcal{V}_n > d$ for all $n \geq m$.  Thus, $\mcal{V}_n > \mcal{V}_m + (n - m) d$ for $n
> m$. It follows that
{\small \bee \sum_{n \in \bb{N} } \exp \li (  \li [ \f{ \varphi(\se) }{ \se } -  \f{ \varphi(s) }{ s } \ri ] \se \mcal{V}_n \ri ) & < & \sum_{n
= 1}^{m-1} \exp \li ( \li [ \f{ \varphi(\se) }{ \se } -  \f{ \varphi(s) }{ s } \ri ] \se \mcal{V}_n \ri ) + \sum_{n = m}^\iy \exp \li ( \li [
\f{ \varphi(\se) }{ \se } -  \f{ \varphi(s) }{ s } \ri ] \se [ \mcal{V}_m + (n - m) d ] \ri )\\
& = &  \sum_{n = 1}^{m-1} \exp \li ( \li [ \f{ \varphi(\se) }{ \se } -  \f{ \varphi(s) }{ s } \ri ] \se \mcal{V}_n \ri ) + \f{ \exp \li ( \li [
\f{ \varphi(\se) }{ \se } -  \f{ \varphi(s) }{ s } \ri ] \se \mcal{V}_m \ri ) }{ 1 - \exp \li ( \li [ \f{ \varphi(\se) }{ \se } -  \f{
\varphi(s) }{ s } \ri ] \se d \ri ) } < \iy. \eee}  This implies that {\small $\sum_{n \in \bb{N} } \Pr \{ Y_n \geq \ga \}$} is finite. It
follows from Borel-Cantelli lemma that $\Pr \{ \cap_{n = 1}^\iy \cup_{k \geq n} [ Y_k \geq \ga ] \} = 0$ and thus $\lim_{n \to \iy } Y_n = 0$
almost surely for $s \in (0, b)$.   This proves the claim that for $s \in (0, b)$, $\lim_{t \to \iy} Y_t = 0$ almost surely. Similarly, we can
show that for $s \in (-a, 0)$, $\lim_{t \to \iy} Y_t = 0$ almost surely.  Since for $s \in (-a, 0) \cup (0, b)$, $(Y_t, \mscr{F}_t)_{t \in
\bb{R}^+}$ is a non-negative continuous martingale which converges to $0$, the proof of the theorem can be completed by applying Theorem
\ref{Them5}.

\section{Proof of Theorem \ref{Them9}} \la{Them9_app}

Define $W_t = \exp ( s (X_t - X_0) - \varphi (s) \mcal{V}_t )$ for $ t \geq 0$ and $s \in (-a, b)$. Then, for all $s \in (-a, b)$ and arbitrary
$t^\prime \geq t \geq 0$, we have {\small \bee &  & \bb{E} [ W_{t^\prime} \mid \mscr{F}_t] = \bb{E} \li [  \exp ( s ( X_{t^\prime} - X_0 ) -
\varphi (s) \mcal{V}_{t^\prime} ) \mid \mscr{F}_t \ri ]
 = \bb{E} \li [  \exp ( s (X_{t^\prime} - X_t) - \varphi (s) (\mcal{V}_{t^\prime} - \mcal{V}_t)  ) \; W_t \mid \mscr{F}_t \ri ]\\
&  & =  W_t  \exp ( - \varphi (s) (\mcal{V}_{t^\prime} - \mcal{V}_t) ) \; \bb{E} \li [  \exp ( s (X_{t^\prime} - X_t)  ) \mid \mscr{F}_t \ri ]
\leq W_t. \eee} Hence, for any $s \in (-a, b)$, $( W_t, \mscr{F}_t )_{t \in \bb{R}^+ }$ is a super-martingale with $\bb{E} [W_0] = \bb{E} [ \exp
( - \varphi (s) \mcal{V}_0 ) ] \leq 1$.  By the assumption on the continuity of the sample paths of $\{ s X_t - \varphi (s) \mcal{V}_t \}_{t >
0}$, we have that almost all sample paths of $( W_t )_{t \in \bb{R}^+}$ is right-continuous.

To prove (\ref{gena12}), note that for any  $s \in (0, a)$ and real number $\ga > 0$,  \bel & & \Pr \li \{
\inf_{t > 0} \li [ X_t - X_0 + \ga \mcal{V}_\tau + \f{ \varphi(-s) }{s} (\mcal{V}_t - \mcal{V}_\tau) \ri ] \leq 0 \ri \} \nonumber\\
&  & = \Pr \li \{ \inf_{t > 0} \li [ X_t - X_0 + \ga \mcal{V}_\tau + \f{ \varphi(-s) }{s} (\mcal{V}_t - \mcal{V}_\tau) \ri ] s \leq 0 \ri \} \nonumber\\
&  & = \Pr \li \{ \inf_{t > 0} \li [ s ( X_t - X_0 ) + \varphi(- s) \mcal{V}_t + \ga s \mcal{V}_\tau  - \varphi(- s) \mcal{V}_\tau \ri ] \leq 0 \ri \} \nonumber\\
&  & = \Pr \li \{ \sup_{t > 0} \li [ - s ( X_t - X_0 )- \varphi(- s) \mcal{V}_t - \ga s \mcal{V}_\tau  + \varphi(- s) \mcal{V}_\tau \ri ] \geq 0 \ri \} \nonumber\\
&  & = \Pr \li \{ \sup_{t > 0} \li [ - s ( X_t - X_0 ) - \varphi(- s) \mcal{V}_t \ri ] \geq \ga s \mcal{V}_\tau  - \varphi(- s) \mcal{V}_\tau  \ri \}\nonumber\\
&  & = \Pr \li \{ \sup_{t > 0} W_t  \geq \exp \li ( \ga s \mcal{V}_\tau  - \varphi(- s) \mcal{V}_\tau \ri ) \ri \} \la{defma896}\\
&  & \leq \exp \li ( \varphi(- s) \mcal{V}_\tau - \ga s \mcal{V}_\tau   \ri ) \la{marineqa8986}\\
&  &  = \li [ \exp \li ( \varphi(- s) - \ga s  \ri ) \ri ]^{\mcal{V}_\tau}. \nonumber \eel Here, we have used the definition of $W_t$ in
(\ref{defma896}). The inequality (\ref{marineqa8986}) follows from the super-martingale inequality.  This proves (\ref{gena12}).

To prove (\ref{genb12}), note that for any $s \in (0, b)$ and real number $\ga > 0$, \bel &  & \Pr \li \{
\sup_{t > 0} \li [ X_t - X_0 - \ga \mcal{V}_\tau - \f{ \varphi(s) }{s} (\mcal{V}_t - \mcal{V}_\tau) \ri ] \geq 0 \ri \} \nonumber\\
&  & = \Pr \li \{ \sup_{t > 0} \li [ X_t - X_0 - \ga \mcal{V}_\tau - \f{ \varphi(s) }{s} (\mcal{V}_t - \mcal{V}_\tau) \ri ] s \geq 0 \ri \} \nonumber\\
&  & = \Pr \li \{ \sup_{t > 0} \li [ s ( X_t - X_0 ) - \varphi(s) \mcal{V}_t - \ga s \mcal{V}_\tau  + \varphi(s) \mcal{V}_\tau \ri ] \geq 0 \ri \} \nonumber\\
&  & = \Pr \li \{ \sup_{t > 0} \li [ s ( X_t - X_0 ) - \varphi(s) \mcal{V}_t \ri ] \geq \ga s \mcal{V}_\tau  - \varphi(s) \mcal{V}_\tau  \ri \}\nonumber\\
&  & = \Pr \li \{ \sup_{t > 0} W_t  \geq \exp \li ( \ga s \mcal{V}_\tau  - \varphi(s) \mcal{V}_\tau \ri ) \ri \} \la{defm336}\\
&  & \leq \exp \li ( \varphi(s) \mcal{V}_\tau - \ga s \mcal{V}_\tau  \ri ) \la{marineq339}\\
&  &  = \li [ \exp \li ( \varphi(s) - \ga s  \ri ) \ri ]^{\mcal{V}_\tau}. \nonumber \eel  Here, we have used the definition of $W_t$ in
(\ref{defm336}). The inequality (\ref{marineq339}) follows from the super-martingale inequality.  This proves (\ref{genb12}).

Before proving (\ref{discovera99812}) and (\ref{discoverb89912}), we shall show (\ref{discovera88312}) and (\ref{discoverb9812}). Note that, as
a consequence of $0 \leq \varphi(s) \leq \ga s$, we have $\varphi(s) ( \mcal{V}_t - \mcal{V}_\tau) \leq \ga s (\mcal{V}_\tau \vee \mcal{V}_t -
\mcal{V}_\tau)$ or equivalently, $\varphi(s) \mcal{V}_t + \ga s \mcal{V}_\tau - \varphi(s) \mcal{V}_\tau \leq \ga s (\mcal{V}_\tau \vee
\mcal{V}_t)$ for any $t > 0$.  This inequality can be written as \be \la{greater} \varphi(s) \mcal{V}_t + s (\eta + \ga \mcal{V}_\tau) -
\varphi(s) \mcal{V}_\tau \leq  \eta s + \ga s (\mcal{V}_\tau \vee \mcal{V}_t). \ee Hence, for any $s \in (0, b)$,
\bel &  & \Pr \li \{ \sup_{ t > 0 } \li [ X_t - X_0 - \eta  - \ga (\mcal{V}_\tau \vee \mcal{V}_t)  \ri ] \geq 0 \ri \} \nonumber\\
&  & = \Pr \li \{ \sup_{ t > 0 } \li [ s ( X_t - X_0 ) - \eta s - \ga s (\mcal{V}_\tau \vee \mcal{V}_t)  \ri ] \geq 0 \ri \} \nonumber\\
&  & \leq \Pr \li \{ \sup_{ t > 0 } \li [ s ( X_t - X_0 ) - \varphi(s) \mcal{V}_t - s (\eta + \ga \mcal{V}_\tau)
 + \varphi(s) \mcal{V}_\tau \ri ] \geq 0 \ri \} \la{great}\\
&  & = \Pr \li \{ \sup_{ t > 0 } \li [ s ( X_t - X_0 ) - \varphi(s) \mcal{V}_t \ri ]
\geq s (\eta + \ga \mcal{V}_\tau)  - \varphi(s) \mcal{V}_\tau  \ri \} \nonumber\\
&  & = \Pr \li \{ \sup_{ t > 0 } W_t  \geq \exp \li ( s (\eta + \ga \mcal{V}_\tau)  - \varphi(s) \mcal{V}_\tau \ri ) \ri \} \la{defm}\\
&  & \leq \exp \li (  \varphi(s) \mcal{V}_\tau - s (\eta + \ga \mcal{V}_\tau) \ri ) \la{marineq}\\
&  &  = e^{- \eta s} \li [ \exp \li (  \varphi(s) - \ga s   \ri ) \ri ]^{\mcal{V}_\tau}. \nonumber \eel    Here, (\ref{great}) follows from
(\ref{greater}). We have used the definition of $W_t$ in (\ref{defm}).  Recall that, for any $s \in (-a, b)$, $( W_t, \mscr{F}_t )_{t \in
\bb{R}^+}$ is a super-martingale with $\bb{E} [W_0] \leq 1$.   The inequality (\ref{marineq}) follows from the super-martingale inequality. This
proves (\ref{discoverb9812}).  The proof of (\ref{discovera88312}) is similar.

Now we are in position to prove (\ref{discovera99812}) and (\ref{discoverb89912}).  In the case that $\{ s \in (0, b): \varphi(s) \leq \ga s \}
$ is empty, (\ref{discoverb89912}) is clearly true, since the infimum is no less than $1$.  In the case that $\{ s \in (0, b): \varphi(s) \leq
\ga s \} $ is not empty, it follows from (\ref{discoverb9812}) that {\small
\[ \Pr \li \{ \sup_{ t > 0 } \li [ X_t - X_0 - \ga (\mcal{V}_\tau \vee \mcal{V}_t) \ri ] \geq 0 \ri \} \leq \inf_{ \{ s \in (0, b): \varphi(s)
\leq \ga s \} }  \li [ \exp \li ( \varphi(s) - \ga s  \ri ) \ri ]^{\mcal{V}_\tau} = \inf_{ s \in (0, b) }  \li [ \exp \li ( \varphi(s) - \ga s
\ri ) \ri ]^{\mcal{V}_\tau}.
\]}
This proves (\ref{discoverb89912}).  The proof of (\ref{discovera99812}) is similar.

To prove (\ref{discoverb9899add}), note that, for any $s \in \mscr{A}$,
\bee &  & \Pr \li \{ \inf_{ t > 0 } \li [ X_t - X_0 + \eta  + \ga  \mcal{V}_t  \ri ] \leq 0 \ri \}
= \Pr \li \{ \sup_{ t > 0 } \li [ (- s) ( X_t - X_0 ) - \eta s - \ga s \mcal{V}_t  \ri ] \geq 0 \ri \} \nonumber\\
&  & \leq  \Pr \li \{ \sup_{ t > 0 } \li [ (- s) ( X_t - X_0 ) - \varphi(- s) \mcal{V}_t \ri ] \geq s \eta   \ri \}  = \Pr \li \{ \sup_{ t > 0 }
W_t \geq \exp  ( s \eta  ) \ri \}  \leq \exp  (  - s \eta  ),  \eee    where $W_t = \exp \li [ (- s) ( X_t - X_0 ) - \varphi(- s) \mcal{V}_t \ri
]$ and  the last inequality follows from the super-martingale inequality.

To prove (\ref{discoverb989912}), note that, for any $s \in \mscr{B}$,
\bee &  & \Pr \li \{ \sup_{ t > 0 } \li [ X_t - X_0 - \eta  - \ga  \mcal{V}_t  \ri ] \geq 0 \ri \}
= \Pr \li \{ \sup_{ t > 0 } \li [ s ( X_t - X_0 ) - \eta s - \ga s \mcal{V}_t  \ri ] \geq 0 \ri \} \nonumber\\
&  & \leq  \Pr \li \{ \sup_{ t > 0 } \li [ s ( X_t - X_0 ) - \varphi(s) \mcal{V}_t \ri ] \geq s \eta   \ri \}  = \Pr \li \{ \sup_{ t > 0 } W_t
\geq \exp  ( s \eta  ) \ri \}  \leq \exp  (  - s \eta  ),  \eee    where $W_t = \exp \li [ s ( X_t - X_0 ) - \varphi(s) \mcal{V}_t \ri ]$ and
the last inequality follows from the super-martingale inequality.

We shall show statement (I). For simplicity of notation, define $g(s) = \varphi(s) - \ga s$ for $s \in (0, b)$. To show (\ref{discoverb12}),
recall the assumption that $\varphi(s)$ is a non-negative, continuous function smaller than $\ga |s|$ at a neighborhood of $0$. Hence,  there
exists a number $\ep \in (0, b)$ such that $\varphi(s) < \ga s$ for $s \in (0, \ep ]$.  Since $g(\ep) < 0 = g(0)$, it must be true that either
the infimum of $g(s)$ is attained at some $s^\star \in (0, b)$ or $\inf_{ s \in (0, b) } g(s) = \lim_{s \uparrow b} g(s) < g(\ep) < 0$. In the
former case, (\ref{discoverb12}) of statement (I) is true as a consequence of (\ref{genb12}). In the latter case, we can define $\varphi(b) =
\lim_{s \uparrow b} \varphi(s)$.  Then, $b \lim_{s \uparrow b} \f{ \varphi(s) }{s} = \varphi(b) < \ga b$ and $\varphi(b) - \ga b = \inf_{s \in
(0, b)} [ \varphi(s) - \ga s ]$. Consider $W_t = \exp ( b (X_t - X_0) - \varphi (b) \mcal{V}_t )$ for $ t \geq 0$. For arbitrary $t^\prime \geq
t \geq 0$, making use of  Fatou's lemma, we have {\small \bee &  & \bb{E} \li [ \exp ( b (X_{t^\prime} - X_t)  ) \mid \mscr{F}_t \ri ] = \bb{E}
\li [ \liminf_{s \uparrow b} \exp ( s (X_{t^\prime} -
X_t) ) \mid \mscr{F}_t \ri ] \leq \liminf_{s \uparrow b} \bb{E} \li [  \exp ( s (X_{t^\prime} - X_t)  ) \mid \mscr{F}_t \ri ]\\
&  & \leq \liminf_{s \uparrow b} \exp ( \varphi(s) (\mcal{V}_{t^\prime} - \mcal{V}_t)  ) = \exp ( \varphi(b) (\mcal{V}_{t^\prime} - \mcal{V}_t)
) \eee}  and consequently,  {\small \bee & & \bb{E} [ W_{t^\prime} \mid \mscr{F}_t] = \bb{E} \li [ \exp ( b ( X_{t^\prime} - X_0 ) - \varphi (b)
\mcal{V}_{t^\prime} )  \mid \mscr{F}_t \ri ]
 = \bb{E} \li [  \exp ( b (X_{t^\prime} - X_t) - \varphi (b) (\mcal{V}_{t^\prime} - \mcal{V}_t)  ) \; W_t \mid \mscr{F}_t \ri ]\\
&  & =  W_t  \exp ( - \varphi (b) (\mcal{V}_{t^\prime} - \mcal{V}_t) ) \; \bb{E} \li [  \exp ( b (X_{t^\prime} - X_t)  ) \mid \mscr{F}_t \ri ]
\leq W_t \eee} almost surely.  This implies that $( W_t, \mscr{F}_t )_{t \in \bb{R}^+}$ is a super-martingale with $\bb{E} [W_0] = \bb{E} [ \exp
( - \varphi (b) \mcal{V}_0 ) ] \leq 1$. So, for any real number $\ga > 0$, \bee &  & \Pr \li \{
\sup_{t > 0} \li [ X_t - X_0 - \ga \mcal{V}_\tau - \lim_{s \uparrow b} \f{ \varphi(s) }{s} (\mcal{V}_t - \mcal{V}_\tau) \ri ] \geq 0 \ri \} \\
&  & = \Pr \li \{ \sup_{t > 0} \li [ X_t - X_0 - \ga \mcal{V}_\tau - \lim_{s \uparrow b} \f{ \varphi(s) }{s} (\mcal{V}_t - \mcal{V}_\tau) \ri ] b \geq 0 \ri \} \\
&  & = \Pr \li \{ \sup_{t > 0} \li [ b ( X_t - X_0 ) - \varphi(b) \mcal{V}_t - \ga b \mcal{V}_\tau  + \varphi(b) \mcal{V}_\tau \ri ] \geq 0 \ri \}\\
&  & = \Pr \li \{ \sup_{t > 0} \li [ b ( X_t - X_0 ) - \varphi(b) \mcal{V}_t \ri ] \geq \ga b \mcal{V}_\tau  - \varphi(b) \mcal{V}_\tau  \ri \}\\
&  & = \Pr \li \{ \sup_{t > 0} W_t  \geq \exp \li ( \ga b \mcal{V}_\tau  - \varphi(b) \mcal{V}_\tau \ri ) \ri \} \leq \li [ \exp \li (
\varphi(b) - \ga b \ri ) \ri ]^{\mcal{V}_\tau} = \inf_{s \in (0, b)} \li [ \exp \li ( \varphi(s) - \ga s \ri ) \ri ]^{\mcal{V}_\tau}. \eee This
establishes (\ref{discoverb12}).  Now we shall show that $0 < \ba(\ga) < \ga$.  Clearly, $\ba(\ga)$ is positive.  In the case that the infimum
of $g(s)$ is attained at some $s^\star \in (0, b)$, we have $g(s^\star) < g (\ep) < 0$, which implies that  $\ba(\ga) < \ga$.  In the case that
$\inf_{ s \in (0, b) } g(s) = \lim_{s \uparrow b} g (s) < 0$, we have $\f{\lim_{s \uparrow b} \varphi(s)}{b} = \lim_{s \uparrow b}
\f{\varphi(s)}{s} < \ga$. So, in both cases, $0 < \ba(\ga) < \ga$.

In a similar manner, we can show (\ref{discovera12}) and  the inequality $0 < \al(\ga) < \ga$.

Finally, we need to show statement (II). This is clearly true for the case that $\lim_{s \uparrow b} \f{ \varphi(s) }{s} \leq \ga$.   It
suffices to consider the case that  $\lim_{s \uparrow b} \f{ \varphi(s) }{s} > \ga$.   In this case, by the assumption that $\f{ \varphi(s)
}{|s|}$ is monotonically increasing with respect to $|s| > 0$,  there exists a unique number $b^* \in (0, b)$ such that $\f{ \varphi(b^*) }{b^*}
= \ga$. Therefore, $\inf_{ \{ s \in (0, b): \varphi(s) \leq \ga s \} } \exp \li ( \varphi(s) - \ga s \ri ) = \inf_{ s \in (0, b^*] } \exp \li (
\varphi(s) - \ga s \ri )$. It follows from (\ref{discoverb9812}) that (\ref{Julya}) is true.  In a similar manner, we can show (\ref{Julyb}).
This concludes the proof of the theorem.

\section{An Upper Bound for the Moment Generating Function of a Uniform Random Variable} \la{app_uniform}

In this appendix, we shall establish the following result.

\beT \la{thmUniform}

Let $Y$ be a random variable uniformly distributed over $[- \f{1}{2}, \f{1}{2}]$.  Then, $\bb{E} [ \exp( s Y) ] \leq \exp \li ( \f{s^2}{24} \ri
)$ for any $s \in \bb{R}$.  That is, the moment generating function of $Y$ is bounded from above by $\exp \li ( \f{s^2}{24} \ri )$. \eeT

\bpf Note that the moment generating function of $Y$ is $\bb{E} [ \exp( s Y) ] = \f{  (e^s - 1) e^{- \f{s}{2}} }{s}$.  We want to show that
$g(s) < \exp \li (  \f{s^2}{24} \ri )$ for any $s \in \bb{R}$.  Define $h(s) = s \li [  \bb{E} [ \exp( s Y) ] -  \exp \li (  \f{s^2}{24} \ri )
\ri ]$.  Then, $h(s) = (e^s - 1) e^{- \f{s}{2}} - s \exp \li (  \f{s^2}{24} \ri )$.  It can be checked that the derivative of $h(s)$ is
$h^\prime (s) = u(s) - v(s)$, where
\[
u(s) = \f{1}{2} \li ( e^{ \f{s}{2}} + e^{- \f{s}{2}}  \ri ), \qqu v(s)  = \li ( 1 + \f{s^2}{12} \ri ) \exp \li (  \f{s^2}{24} \ri ).
\]
Using Taylor series expansion formula, we can write
\[
u(s) = 1 + \sum_{i = 1}^\iy \f{ s^{2 i}  } { 4^{i} (2 i)! }, \qqu  v(s) = 1 +  \sum_{i = 1}^\iy \li [  \f{ \li (  \f{s^2}{24} \ri )^i  }{i!} +
 \f{ \li (  \f{s^2}{24} \ri )^{i-1} \f{s^2}{12} }{ (i -1)!} \ri ] = 1 +
\sum_{i = 1}^\iy  \f{ 1 + 2 i} { 24^i \; i!} s^{2 i}.
\]
Since $h^\prime (0) = u(0) - v(0)$, to show the theorem, it suffices to show that $u(s) < v(s)$.  This can be accomplished by proving $\f{ s^{2
i}  } { 4^{i} (2 i)! } < \f{ 1 + 2 i} { 24^i \; i!} s^{2 i}$ for $i = 1, 2, \cd$, or equivalently,
\[
\f{1}{ (2 i)! } \leq \f{ 1 + 2 i} {6^i i!}, \qqu i = 1, 2, \cd.
\]
To do so, define the ratio of $\f{ 1 + 2 i} {6^i i!}$ to $\f{1 }{ (2 i)! }$ as $f(i)$.  It can be checked that $f (i) = \f{ (2 i)! ( 1 + 2 i) }
{ 6^i i!  }$ for $i = 1, 2, \cd$.  Since $f(1) \geq 1$ and $\f{ f(i + 1)  }{ f(i) } = 1 + \f{2 i}{3} > 1$ for $k \geq 1$, we have $f (i)
> 1 $ for all $i > 1$.  This implies that $u(s) < v(s)$ for $s \in \bb{R}$.  The proof of the theorem is thus completed.

\epf

\section{Proof of Corollary \ref{expmax}} \la{expmax_app}

Let $\psi (.)$ be the inverse function of $u (.)$ such that $u (\psi(\ze)) = \ze$ for $\ze \in \{ u (\se): \se \in \Se \}$.  Define $h(\ze) = v
(\psi(\ze))$ for $\ze \in \{ u (\se): \se \in \Se \}$.   Putting $\zeta = u (\se)$, we have $\bb{E} \li [ \exp ( s Y  ) \ri ] = \exp \li ( h
(\zeta + s) - h (\zeta) \ri )$.  Define $\varphi(s) = h (\zeta + s) - h (\zeta) - \se s$ and \[ Z_n = X_n - n \se, \qqu \mcal{V}_n = n
\]
for $n \in \bb{N}$.  For $n \in \bb{N}$, let $\mscr{F}_n$ denote the $\si$-algebra generated by $X_1, \cd,  X_n$. Then,
\[
\bb{E} [ \exp ( s (Z_n  - Z_m) ) \mid \mscr{F}_m ] = \exp ( (\mcal{V}_n - \mcal{V}_m) \varphi(s) )
\]
for $m, n \in \bb{N}$ with $m \leq n$.   Assume that $\se + \ga \in \Se$.  Noting that \be \la{inevip}  \f{ d h(\ze) }{ d \ze } =  \f{ d v }{ d
\psi } \f{ d \psi }{ d \ze } = \psi  \f{ d u }{ d \psi } \f{ d \psi }{ d \ze } =  \psi  \f{ d u }{ d \ze } = \psi (\ze), \ee we have
\[
\f{ d [ \varphi(s) - \ga s] }{ d s } = \f{ d h (\zeta + s)  }{ d s } - (\se + \ga) =  \psi (\zeta + s) - (\se + \ga).
\]
Define $s^* = u (\se + \ga) - u (\se)$. Invoking the definition that $\ze = u(\se)$, we have $\zeta + s^* = u (\se + \ga)$, which implies that
$\psi (\zeta + s^*) = \se + \ga$.  Therefore,
\[
\li. \f{ d [ \varphi(s) - \ga s] }{ d s } \ri |_{s = s^*} = \psi (\zeta + s^*) - (\se + \ga) = 0.
\]
Since $\varphi(s) - \ga s$ is a convex function of $s$, the infimum of $\varphi(s) - \ga s$ is attained at $s = s^*$.  Note that
\[
\varphi(s^*) = h (\zeta + u (\se + \ga) - u (\se) ) - h (\zeta) - \se [ u (\se + \ga) - u (\se) ] = v(\se + \ga) - v(\se) - \se [ u (\se + \ga)
- u (\se) ].
\]
Thus, \bel &  & \varphi( s^* ) - \ga s^* =  v(\se + \ga) - v(\se) + (\se + \ga) [  u (\se) - u (\se + \ga) ], \la{use33}\\
&  & \f{\varphi(s^*)}{s^*} =  \f{ v(\se + \ga) - v(\se) }{ u (\se + \ga) - u (\se) } - \se \in (0, \ga). \la{use88} \eel Applying Theorem
\ref{Them9}, we have {\small \be \la{use99} \Pr \li \{ \sup_{n > 0} \li [ X_n - n \se  - \ga(n \vee m) \ri ] \geq  0 \ri \}  \leq \Pr \li \{
\sup_{n > 0} \li [ X_n - n \se - m \ga - \f{\varphi(s^*)}{s^*} (n - m) \ri ] \geq 0 \ri \} \leq \exp ( m [  \varphi( s^* ) - \ga s^*  ] ).  \ee}
Making use of (\ref{use33}), (\ref{use88}), (\ref{use99}) and the definitions of $\rho$ and $\mscr{M}$, we have that $\Pr \{ \sup_{n > 0} [ X_n
- n \se  - \ga(n \vee m) ] \geq  0 \} \leq \Pr \{ \sup_{n > 0} [ X_n -  \ro (\se + \ga, \se, m, n) ] \geq  0 \} \leq  [ \mscr{M} (\se + \ga, \se
) ]^m$ provided that  $\se + \ga \in \Se$.  This proves (\ref{ineqgood8}).

Now we shall show (\ref{ineqgood9}).  Assume that $\se - \ga \in \Se$.  By virtue of (\ref{inevip}), we have
\[
\f{ d [ \varphi( - s) - \ga s] }{ d s } = \f{ d h (\zeta - s)  }{ d s} + (\se - \ga) =  - \psi (\zeta - s) + (\se - \ga).
\]
Define $s^\star = u (\se) - u (\se - \ga)$. Invoking the definition that $\ze = u(\se)$, we have $\zeta - s^\star = u (\se - \ga)$, which
implies that $\psi (\zeta - s^\star) = \se - \ga$.  Therefore,
\[
\li. \f{ d [ \varphi(-s) - \ga s] }{ d s } \ri |_{s = s^\star} = \se - \ga - \psi (\zeta - s^\star) = 0.
\]
Since $\varphi(- s) - \ga s$ is a convex function of $s$, the infimum of $\varphi(-s) - \ga s$ is attained at $s = s^\star$.  Note that
\[
\varphi(- s^\star) = h (\zeta - u (\se) + u (\se - \ga) ) - h (\zeta) + \se [ u (\se) - u (\se - \ga) ] = v(\se - \ga) - v(\se) + \se [ u (\se)
- u (\se - \ga) ].
\]
Thus, \bel &  & \varphi( - s^\star ) - \ga s^\star =  v(\se - \ga) - v(\se) + (\se - \ga) [  u (\se) - u (\se - \ga) ], \la{use33b}\\
&  & \f{\varphi(- s^\star)}{s^\star} =  \f{ v(\se - \ga) - v(\se) }{ u (\se ) - u (\se - \ga) } + \se \in (0, \ga). \la{use88b} \eel Applying
Theorem \ref{Them9}, we have {\small \be \la{use99b} \Pr \li \{ \inf_{n > 0} \li [ X_n - n \se  + \ga(n \vee m) \ri ] \leq  0 \ri \}  \leq \Pr
\li \{ \inf_{n > 0} \li [ X_n - n \se + m \ga + \f{\varphi(- s^\star)}{s^\star} (n - m) \ri ] \leq 0 \ri \} \leq \exp ( m [  \varphi( - s^\star
) - \ga s^\star ] ).  \ee} Making use of (\ref{use33b}), (\ref{use88b}), (\ref{use99b}) and the definitions of $\rho$ and $\mscr{M}$, we have
that $\Pr \{ \inf_{n > 0} [ X_n - n \se  + \ga(n \vee m) ] \leq  0 \} \leq \Pr \{ \inf_{n > 0} [ X_n -  \ro (\se - \ga, \se, m, n) ] \leq  0 \}
\leq  [ \mscr{M} (\se - \ga, \se ) ]^m$ provided that  $\se - \ga \in \Se$.  This proves (\ref{ineqgood9}).  The proof of the theorem is thus
completed.

\section{Proof of Corollary \ref{CorAzuma}} \la{CorAzuma_app}

For simplicity of notations, define $A_{t, \tau} = Y_\tau - X_\tau - \sq{\mcal{V}_t - \mcal{V}_\tau}$ and $B_{t, \tau} = Y_\tau - X_\tau +
\sq{\mcal{V}_t - \mcal{V}_\tau} + \de$, where $\de$ is a positive number introduced for the purpose of ensuring $B_{t, \tau} - A_{t, \tau} > 0$.
By the assumption that $|X_t - Y_\tau|^2 \leq \mcal{V}_t - \mcal{V}_\tau$, we have that $A_{t, \tau} \leq X_t - X_\tau < B_{t, \tau}$ almost
surely.  For $t \geq \tau \geq 0$, we have {\small \bel \bb{E} \li [  e^{s (X_t - X_\tau) } \mid \mscr{F}_\tau \ri ] & \leq & \bb{E} \li [ \f{
B_{t, \tau} - (X_t - X_\tau)} {B_{t, \tau} - A_{t, \tau}} e^{s A_{t, \tau}} +
\f{(X_t - X_\tau) - A_{t, \tau}} {B_{t, \tau} - A_{t, \tau}} e^{s B_{t, \tau}} \mid \mscr{F}_\tau \ri ] \la{linearconvex} \\
& = & \f{ B_{t, \tau} \; e^{s A_{t, \tau}}  - A_{t, \tau} \; e^{s B_{t, \tau}} } {B_{t, \tau} - A_{t, \tau}} + \f{e^{s B_{t, \tau}} - e^{s A_{t,
\tau}}} {B_{t, \tau} - A_{t, \tau}} \bb{E} \li [ X_t - X_\tau  \mid
\mscr{F}_\tau \ri ] \la{conmeasure}\\
& \leq & \f{ B_{t, \tau} \; e^{s A_{t, \tau}}  - A_{t, \tau} \; e^{s B_{t, \tau}} } {B_{t, \tau} - A_{t, \tau}} \la{consupm}\\
& \leq & \exp \li ( \f{ s^2 (B_{t, \tau} - A_{t, \tau})^2 }{8} \ri ) \la{ineH} \\
&  = & \exp \li ( \f{ s^2 (\sq{\mcal{V}_t - \mcal{V}_\tau} + \f{\de}{2})^2 }{2} \ri ) \la{desure} \eel} almost surely, where
(\ref{linearconvex}) follows from the convexity of the exponential function,  (\ref{conmeasure}) follows from the fact that $A_{t, \tau}$ and
$B_{t, \tau}$ are measurable in $\mscr{F}_\tau$, (\ref{consupm}) follows from the assumption that $(X_t, \mscr{F}_t)_{t \in \bb{R}^+ }$ is a
super-martingale, (\ref{ineH}) follows from the inequality $y e^x - x e^y \leq (y - x) \exp ( |y - x|^2 \sh 8)$ for $y \geq x$.  Since
(\ref{desure}) holds almost surely for any $\de > 0$, it must be true that for any $s \in \bb{R}$,
\[
\bb{E} \li [  e^{s (X_t - X_\tau) } \mid \mscr{F}_\tau \ri ] \leq \lim_{\de \downarrow 0} \exp \li ( \f{ s^2 (\sq{\mcal{V}_t - \mcal{V}_\tau} +
\f{\de}{2})^2 }{2} \ri )  = \exp \li ( \f{ (\mcal{V}_t - \mcal{V}_\tau) s^2}{2} \ri )
\]
almost surely.  Therefore, we have shown that $\bb{E} \li [  e^{s (X_t - X_\tau) } \mid \mscr{F}_\tau \ri ] \leq \exp \li ( (\mcal{V}_t -
\mcal{V}_\tau) \varphi(s) \ri )$ with $\varphi(s) = \f{s^2}{2}$.  Clearly, $\varphi(0) = 0$ and $\varphi(s)$ is convex. Moreover, for $\vep >
0$, $\inf_{s \in (0, \iy)} [ \varphi(s) - \vep s ]$ is attained at $s = s^*$ with $s^* = \vep$.   It can be checked that $\varphi(s^*) - \vep
s^* = - \f{\vep^2}{2}$ and $\f{ \varphi(s^*) }{s^*} = \f{\vep}{2}$.  Hence, applying Theorem \ref{Them9}, we have
\[
\Pr \li \{   \sup_{t > 0}  \li [  X_t - X_0 - \vep \mcal{V}_\tau - \f{\vep}{2} \li ( \mcal{V}_t - \mcal{V}_\tau \ri ) \ri ] \geq 0 \ri \} \leq
\exp \li ( - \f{ \vep^2 \mcal{V}_\tau } { 2 } \ri ).
\]
Substituting $\ga$ in the above inequality with $\vep \mcal{V}_\tau$ yields (\ref{hoe9a}).  In a similar manner, we can show (\ref{hoe9b}).

To show (\ref{hoe9c}), note that \bee &  &  \li \{   \sup_{t > 0} \li [  X_t - X_0 - \f{\ga}{2} \li ( 1 + \f{\mcal{V}_t}{\mcal{V}_\tau} \ri )
\ri ] < 0 \ri \} \bigcap \li \{   \inf_{t > 0} \li [  X_t - X_0 + \f{\ga}{2} \li ( 1 + \f{\mcal{V}_t}{\mcal{V}_\tau} \ri ) \ri ] < 0   \ri
\}\\
&   & \subseteq \li \{   \sup_{t > 0} \li [  \li | X_t - X_0 \ri | - \f{\ga}{2} \li ( 1 + \f{\mcal{V}_t}{\mcal{V}_\tau} \ri ) \ri ] < 0 \ri \}
\eee Making use of this observation, Bonferroni's inequality,  inequalities (\ref{hoe9a}) and  (\ref{hoe9b}), we have
\[
\Pr \li \{   \sup_{t > 0} \li [  \li | X_t - X_0 \ri | - \f{\ga}{2} \li ( 1 + \f{\mcal{V}_t}{\mcal{V}_\tau} \ri ) \ri ] < 0 \ri \} \geq 1 - 2
\exp \li ( - \f{\ga^2} { 2 \mcal{V}_\tau } \ri ),
\]
from which (\ref{hoe9c}) immediately follows.   This completes the proof of the theorem.

\section{Proof of Corollary \ref{CBB}} \la{CBB_app}

Define  $g (s) = \sum_{k = 2}^\iy \f{s^{k - 2}}{k!}$.  Then, $g(s) = \f{e^s - 1 - s}{s^2}$ for $s \neq 0$ and $g(0) = \f{1}{2}$.    It can be
shown that $g(s)$ is increasing with respect to $s$.  It is well known that $g (s)  \leq   \f{1}{2 ( 1 - s \sh 3)}$ for $s \in (0, 3)$.

\subsection{ Proof of (\ref{beneet})  }

Making use of the assumptions of the theorem and following the techniques of \cite{Bennett, Fan}, we have that {\small \bee & & \bb{E} \li [
e^{s (X_n - X_{n - 1} - a_n)} \mid \mscr{F}_{n - 1} \ri ] = \bb{E}
\li [  \sum_{k = 0}^\iy \f{s^k}{k!} (X_n - X_{n - 1} - a_n)^k  \mid \mscr{F}_{n - 1} \ri ] \\
&  & = 1 - s a_n +  \bb{E} \li [  \sum_{k = 2}^\iy \f{s^k}{k!} (X_n - X_{n - 1} - a_n)^k  \mid \mscr{F}_{n - 1} \ri ]
\leq   1 - s a_n + s^2 g (s b) \; \bb{E} \li [   (X_n - X_{n - 1} - a_n)^2   \mid \mscr{F}_{n - 1} \ri ]\\
&  & = 1 - s a_n + s^2 g (s b) \li ( \bb{E} \li [   (X_n - \bb{E} [X_n \mid \mscr{F}_{n - 1}] )^2    \mid \mscr{F}_{n - 1} \ri ] + a_n^2
\ri )  =  1 - s a_n + s^2 g (s b) \li [ \mrm{Var} ( X_n \mid \mscr{F}_{n - 1} ) + a_n^2 \ri ] \qqu \\
&  &  \leq  1 - s a_n + s^2 g (s b)  \; ( \si_n^2 + a_n^2 ) \leq \exp \li ( - s a_n + s^2 g (s b)  \; ( \si_n^2 + a_n^2 ) \ri ) \eee} almost
surely.  This implies that \be \la{vip988}
 \bb{E} [ e^{s  (X_n - X_{n - 1}) } \mid \mscr{F}_{n - 1} ]  \leq \exp \li ( s^2 g(s b) (\mcal{V}_n -
\mcal{V}_{n - 1}) \ri ) = \exp \li (  \f{ (e^{s b} - 1 - s b)}{b^2} (\mcal{V}_n - \mcal{V}_{n - 1}) \ri ) \ee almost surely.  Thus, we have that
$\bb{E} [ e^{s  (X_n - X_{n - 1}) } \mid \mscr{F}_{n - 1} ] \leq \exp \li ( (\mcal{V}_n - \mcal{V}_{n - 1}) \varphi(s) \ri )$ almost surely,
where $\varphi(s) = \f{ e^{s b} - 1 - s b }{b^2}$.  Now consider $(X_n)_{n \in \bb{N} }$ in the context of Theorem \ref{Them9}.  Clearly,
$\varphi(0) = \varphi^\prime(0) = 0$ and $\varphi(s)$ is convex.  Let $s^* = \f{ \ln ( 1 + b \ga ) }{b}$.  By differentiation, it can be shown
that $\inf_{s \in (0, \iy) } [ \varphi(s) - \ga s]$ is attained at $s = s^*$.   It can be checked that $\varphi (s^*) - \ga s^* = \f{\ga}{b} -
\f{(1 + b \ga) \ln (1 + b \ga) }{b^2}$ and $\f{ \varphi (s^*)  }{s^*} = \f{\ga} {\ln (1 + b \ga)} - \f{1}{b}$. Therefore, applying Theorem
\ref{Them9}, we have  {\small \[ \Pr \li \{ \sup_{n > 0} \li [ X_n - X_0 - \ga \mcal{V}_m  - \li ( \f{\ga} {\ln (1 + b \ga)} - \f{1}{b} \ri )
(\mcal{V}_n - \mcal{V}_m) \ri ] \geq 0 \ri \} \leq \li [ \exp \li ( \f{\ga}{b} - \f{(1 + b \ga) \ln (1 + b \ga) }{b^2} \ri ) \ri ]^{\mcal{V}_m}.
\]}
This completes the proof of (\ref{beneet}).

\subsection{ Proof of (\ref{Bert})  }

Since $g (s)  \leq   \f{1}{2 ( 1 - s \sh 3)}$, it follows from (\ref{vip988})  that
\[
\bb{E} [ e^{s (X_n - X_{n - 1}) } \mid \mscr{F}_{n - 1} ]  \leq \exp \li ( s^2 g(s b) (\mcal{V}_n - \mcal{V}_{n - 1}) \ri )  \leq \exp \li (
\f{s^2}{2 (1 - b s \sh 3)} (\mcal{V}_n - \mcal{V}_{n - 1}) \ri )
\]
almost surely for $s \in (0, \f{3}{b})$.  This implies that for $s \in (0, \f{3}{b})$, $\bb{E} [ e^{s (X_n - X_{n - 1}) } \mid \mscr{F}_{n - 1}
]  \leq \exp \li ( \varphi(s) (\mcal{V}_n - \mcal{V}_{n - 1}) \ri )$ almost surely, where $\varphi(s) = \f{s^2}{2 (1 - b s \sh 3)} $. Applying
 inequality (\ref{genb12}) of  Theorem \ref{Them9} with  $s = \f{1}{ \f{1}{\ga} + \f{b}{3} }$ leads to inequality (\ref{Bert}).

\subsection{ Proof of (\ref{chert8})  }

 Applying inequality (\ref{vip988}) with $a_i = 0$ and $b = 1$, we have that
 \[
\bb{E} [ e^{ s (X_n - X_{n - 1}) } \mid \mscr{F}_{n - 1} ] = \bb{E} [ e^{ s (X_n - X_{n - 1})} ] \leq \exp ( (e^s - 1 - s) \si_n^2) \leq \exp (
s^2 \si_n^2)
 \]
almost surely for $s \in (0, \f{7}{4}]$. Here, we used the fact that $e^s - 1 - s \leq s^2$ for $s \in (0, \f{7}{4}]$.   To see this, consider
function $f(s) \DEF e^s - 1 - s - s^2$.  Note that $f^\prime (s) = e^s - 1 - 2 s$ and $f^{\prime \prime} (s) = e^s - 2$.  Clearly, $f^\prime
(s)$ decreases from $0$ to its minimum at $s = \ln 2$ and then increases for $s \in (\ln 2, \iy)$.  Since $f^\prime (0) = 0$, there exists a
positive number $\ro
> 0$ such that $f^\prime (s)$ is negative for $s \in (0, \ro)$  and positive for $s \in (\ro , \iy)$. Since $f(0) = 0$, it must be true that
$f(s)$ decreases from $0$ to its minimum as $s$ increases from $0$ to $\rho$ and then increases for $s \in (\ro, \iy)$. It follows that if $f(t)
< 0$ for some $t > 0$, then $f(s) < 0 $ for all $s \in (0, t)$.  It can be checked that $f(\f{7}{4}) < 0$, thus, we have $f(s) < 0$ for all  $s
\in (0, \f{7}{4}]$.

So, we have shown that for $s \in (0, \f{7}{4}]$, $\bb{E} [ e^{ s (X_n - X_{n - 1}) } \mid \mscr{F}_{n - 1} ] \leq \exp \li ( (\mcal{V}_n -
\mcal{V}_{n - 1}) \varphi(s) \ri )$ almost surely, where $\varphi(s) = s^2$.  Consider $(X_n)_{n \in \bb{N} }$ in the context of Theorem
\ref{Them9}.  Let $\vep > 0$.  Note that $\f{d \varphi (s) }{d s }= \vep$ if $s = \f{\vep}{2}$. Accordingly, $\ba(\vep) = \f{\vep}{2}$ and $s
(\ba(\vep) - \vep) = - \f{\vep^2}{4}$.   Thus, applying Theorem \ref{Them9},  we have
\[
\Pr \li \{ \sup_{n > 0} \li [ X_n - X_0  - \vep \mcal{V}_m - \f{\vep}{2} (\mcal{V}_n - \mcal{V}_m ) \ri ] \geq 0 \ri \} \leq \exp \li ( -
\f{\vep^2}{4} \ri )^{\mcal{V}_m}
\]
for $\vep \leq \f{7}{2}$.  Letting $\vep \mcal{V}_m = \ga$ in the above inequality leads to (\ref{chert8}).   This completes the proof of the
corollary.

\section{Proof of Corollary \ref{Posmax}} \la{Posmax_app}

Define $Y_t = X_t - \lm t$ for $t \in \bb{R}^+$.    Let $\mcal{V}_t = t$ for $t \in \bb{R}^+$ and let $\varphi(s) = \lm (e^s - 1 - s)$ for $s
\in \bb{R}$. Consider the process  $(Y_t)_{t \in \bb{R}^+}$ in the context of Theorem \ref{Them9}.  Clearly, $(Y_t)_{t \in \bb{R}^+}$ is
actually a right-continuous stochastic process such that $\bb{E} [ \exp ( s (Y_{t^\prime} - Y_t) ) \mid \mscr{F}_t ] \leq \exp (
(\mcal{V}_{t^\prime} - \mcal{V}_t) \varphi(s) )$ almost surely for arbitrary $t^\prime \geq t \geq 0$ and $s \in (-\iy, \iy)$.  Since the
derivative of $\varphi(s) - \ga s$ with respective to $s$ is $\lm e^s - \lm - \ga$, it follows that $\inf_{s \in (0, \iy)} [ \varphi(s) - \ga
s]$ is attained at $s^* = \ln \f{\lm + \ga}{\lm}$. Consequently,
\[
\inf_{s \in (0, \iy)} [ \varphi(s) - \ga s] = \varphi(s^*) - \ga s^* =  \ga - (\lm + \ga) \ln \f{\lm + \ga}{\lm} \qu \tx{and} \qu  \ba(\ga) =
\f{ \varphi(s^*) }{s^*} = \f{ \ga  }{ \ln \f{\lm + \ga}{\lm} } - \lm.
\]
Applying Theorem \ref{Them9} gives
\[
\Pr \li \{ \sup_{t > 0} \li [ Y_t - \ga \tau - \li ( \f{ \ga }{ \ln \f{\lm + \ga}{\lm} } - \lm \ri )  (t - \tau) \ri ] \geq 0 \ri \} \leq \li [
\exp \li ( \ga - (\lm + \ga) \ln \f{\lm + \ga}{\lm} \ri ) \ri ]^\tau,
\]
which implies that {\small $\Pr \{ \sup_{t > 0} [ X_t - (\lm + \ga) \tau - \f{\ga (t - \tau) }{ \ln (1 + \f{\ga}{\lm } ) }  ] \geq 0 \} \leq [ (
\f{\lm}{\lm + \ga} )^{\lm + \ga} e^\ga ]^\tau$} for any $\tau > 0$ and $\ga
> 0$.

On the other hand, since the derivative of $\varphi(- s) - \ga s$ with respective to $s$ is $- \lm e^{-s} + \lm - \ga$, it follows that $\inf_{s
\in (0, \iy)} [ \varphi(- s) - \ga s]$ is attained at $s^\star = \ln \f{\lm}{\lm - \ga}$.  Consequently,
\[
\inf_{s \in (0, \iy)} [ \varphi(- s) - \ga s] =  - \ga - (\lm - \ga) \ln \f{\lm - \ga}{\lm} \qu \tx{and} \qu \al(\ga) = \f{ \varphi(- s^\star)
}{s^\star} = \f{ \ga }{ \ln \f{\lm - \ga}{\lm} } + \lm.
\]
Applying Theorem \ref{Them9} yields
\[
\Pr \li \{ \inf_{t > 0} \li [ Y_t + \ga \tau + \li ( \f{ \ga }{ \ln \f{\lm - \ga}{\lm} } + \lm \ri )  (t - \tau) \ri ] \leq 0 \ri \} \leq \li [
\exp \li ( - \ga - (\lm - \ga) \ln \f{\lm - \ga}{\lm} \ri ) \ri ]^\tau,
\]
which implies that {\small $\Pr \{  \inf_{t
> 0} [ X_t - (\lm - \ga) \tau + \f{ \ga (t - \tau) }{ \ln (1 - \f{\ga}{\lm } ) }   ] \leq 0 \} \leq [ ( \f{\lm}{ \lm - \ga } )^{\lm
- \ga} e^{-\ga} ]^\tau$} for any $\tau > 0$ and $\ga > 0$.  This completes the proof of the corollary.

\end{document}